\documentclass[11pt]{amsart}

\usepackage{amsmath,amsthm,amssymb,amsfonts}
\usepackage{graphicx,color,epsfig,subfigure}
\usepackage{url}
\usepackage{verbatim}
\usepackage{placeins}

\newtheorem{thm}{Theorem}
\newtheorem{defn}[thm]{Definition}
\newtheorem{prop}[thm]{Proposition}
\newtheorem{cor}[thm]{Corollary}
\newtheorem{lem}[thm]{Lemma}

\newtheorem{fact}[thm]{Fact}

\newtheorem{rem}[thm]{Remark}

\newcommand{\vw}{{\rm VW}}
\newcommand{\vb}{{\rm VB}}
\newcommand{\ew}{{\rm EW}}
\newcommand{\eb}{{\rm EB}}

\newcommand{\ews}{\left|{\rm EW}\right|}

\newcommand{\egk}{{\rm EG}(K)}
\newcommand{\xb}{{\mathbf{x}}}
\newcommand{\x}{{\bf x}}
\newcommand{\E}{{\bf E}}

\def\Kt{{K_{2,t}}}
\def\T{{B_{t-2}}}

\def\Hmc{{\mathcal{H}}}

\def\Kmc{{\mathcal{K}}}

\def\dist{{\rm dist}}
\def\forb{{\rm Forb}}
\def\textdef{\textbf}

\def\ed{{\textit{ed}}}

\title{On the edit distance from $K_{2,t}$-free graphs\\ (extended version)}
\author{Ryan R. Martin}
\address{Department of Mathematics, Iowa State University, Ames, Iowa 50011}
\email{rymartin@iastate.edu}
\thanks{The first author's research partially supported by NSF grant DMS-0901008 and by an Iowa State University Faculty Professional Development grant.}
\author{Tracy McKay}
\address{Department of Mathematics, Iowa State University, Ames, Iowa 50011}
\email{tmckay16@iastate.edu}
\thanks{The second author's research partially supported from NSF grant DMS-0901008 and by an Iowa State University Mathematics Department Wolfe Research Fellowship.}

\begin{document}

\maketitle
\begin{abstract}
The edit distance between two graphs on the same vertex set is defined to be the size of the symmetric difference
of their edge sets.  The edit distance function of a hereditary property, $\mathcal{H}$, is a function of $p$,
and measures, asymptotically, the furthest graph of edge density $p$ from $\mathcal{H}$ under this metric. In
this paper, we address the hereditary property $\forb(K_{2,t})$, the property of having no induced copy of the
complete bipartite graph with 2 vertices in one class and $t$ in the other. Employing an assortment of techniques
and colored regularity graph constructions, we are able to determine the edit distance function over the entire
domain $p\in [0,1]$ when $t=3,4$ and extend the interval over which the edit distance function for
$\forb(K_{2,t})$ is known for all values of $t$, determining its maximum value for all odd $t$. We also prove
that the function for odd $t$ has a nontrivial interval on which it achieves its maximum. These are the only
known principal hereditary properties for which this occurs.

In the process of studying this class of functions, we encounter some surprising connections to extremal graph
theory problems, such as strongly regular graphs and the problem of Zarankiewicz.  

This is an extended version of a paper with the same name submitted to the Journal of Graph Theory \cite{jgt_version}.  In particular, this version contains Appendix A, which has tables and graphs pertaining to the hereditary property $\forb(K_{2,t})$ for small $t$, and Appendix B, which has the proofs of Lemma 31, Proposition 32, Proposition 33, and Lemma 34.
\end{abstract}

\section{Introduction}

The study of the edit distance in graphs initially appeared in a paper by Axenovich, K\'ezdy and the first author
\cite{akr2008} and, independently, by Alon and Stav \cite{as2008}. It has several potential applications, such as
to biological consensus trees~\cite{akr2008} and property testing problems in theoretical computer
science~\cite{as2008}.  More recently, interest has been shown in determining the value of the \emph{edit
distance function}, introduced in \cite{BaloghMartin08} by Balogh and the first author. Strategies for
determining this function appear in \cite{MT}, by Marchant and Thomason, as well as in \cite{rm2009}.

Given a hereditary property (that is, a set of graphs closed under vertex deletion and isomorphism), what is the
least number of edge additions or deletions sufficient to make any graph on $n$ vertices  a member of the
property?  What is the behavior of this value as $n\rightarrow\infty$?  In \cite{akr2008}, the binary chromatic
number of a graph $H$  is used as a means of bounding the maximum number of edge additions and deletions (edits)
sufficient to ensure that every $n$-vertex graph has no induced copy of  a single graph $H$.  A hereditary
property which consists of the graphs with no induced copy of $H$ is denoted $\forb(H)$ and is called a
\textdef{principal hereditary property}. The methods in \cite{akr2008} give an asymptotically exact result in
some cases, most notably when $H$ is self-complementary.

Let $G(n,p)$ denote the random graph on $n$ vertices with edge probability $p$.  In \cite{as2008}, a version of
Szemer\'{e}di's regularity lemma is applied to show that, as $n\rightarrow\infty$, the number of edits necessary
to make $G(n,p)$ a member of a given hereditary property approaches, with high probability, the maximum possible
number over all $n$-vertex graphs within $o(n^2)$, so long as $p$ is chosen correctly with respect to the given
hereditary property, $\mathcal{H}$. In fact, the maximum number of edits sufficient to change a density-$p$,
$n$-vertex graph into a member of $\mathcal{H}$ is asymptotically the same as the expected number of edits
required to put $G(n,p)$ into $\mathcal{H}$.

The edit distance function, $\ed_\mathcal{H}(p)$, from \cite{BaloghMartin08}, computes the limit of the maximum normalized edit distance of a density-$p$, $n$-vertex graph from $\Hmc$ as $n\rightarrow\infty$ for all probabilities $p$.  See Definition~\ref{def:distfunc} below.  Not surprisingly, the maximum value of this function occurs at the same $p$ value described in \cite{as2008}.

Marchant and Thomason also explore the value of $1-\ed_\Hmc(p)$ for various hereditary properties in \cite{MT},
developing some insightful results for determining the value of the function in general.  One discovery from
\cite{MT} of particular interest is a relationship between the problem of determining the edit distance function
for $\forb(K_{3,3})$ and constructions by Brown in \cite{Brown} for $K_{3,3}$-free graphs, related to the
Zarankiewicz problem.

More generally, edit distance has been discussed as a potential metric for graph limits, and also as a parameter
for property testing techniques to which graph limits have been applied (see, for example, Borgs, et al.~\cite{BCLSSV2006}).

In this paper, we explore what can be said about the edit distance function for the hereditary property
$\forb(K_{2,t})$, the set of all graphs that do not contain a complete bipartite graph with cocliques of $2$ and
$t$ vertices as an induced subgraph.  In particular, we
\begin{itemize}
   \item Compute the entire edit distance functions for the properties $\forb(K_{2,3})$ and $\forb(K_{2,4})$.
   \item Explore constructions that arise from strongly regular graphs and provide good upper bounds for
   $\ed_{\forb(K_{2,t})}(p)$.  One in particular, derived from the $15$-vertex generalized quadrangle $GQ(2,2)$,
   defines $\ed_{\forb(K_{2,4})}(p)$ for $p\in (1/5,1/3)$.
   \item Compute the edit distance function for the property $\forb(K_{2,t})$ for $p\in [\frac{2}{t+1},1]$.
   \item Show that, for odd $t$, $\ed_{\forb(K_{2,t})}(p)=\frac{1}{t+1}$ for $p\in\left[\frac{2t-1}{t(t+1)},\frac{2}{t+1}\right]$.  These are the only known principal hereditary properties for which the maximum of the edit distance function is achieved on a nontrivial interval.
   \item  Examine the relationship between constructions by F\"uredi~\cite{FurediZar} related to the
   Zarankiewicz problem and the trivial bound $\ed_{\forb(K_{2,t})}(p)\leq p(1-p)$, which achieves the value of the
   function for small values of $p$ when $t=3,4$.  When $t\geq 9$ the constructions by
   F\"uredi~\cite{FurediZar} improve upon this bound for small values of $p$.
   \item Look at constructions derived from powers of cycles that give a general upper bound for some values of $p$ and $t$.
   \item Derive a lower bound for $\ed_{\forb(K_{2,t})}(p)$ that is nontrivial and is achieved for some value of $p$
   if a specified strongly regular graph exists.
   \item Summarize the known bounds for $\ed_{\forb(K_{2,t})}(p)$ for $5\leq t\leq 8$.
\end{itemize}

Prior results and notation come primarily from \cite{rm2009}, as well as previous work: \cite{as2008},
\cite{akr2008}, \cite{BaloghMartin08}, \cite{MT}; however, there are a number of other excellent resources on
related topics. For a more extensive review of this literature, the reader may wish to consult Thomason
\cite{ThomSurv}. We now introduce some important definitions and state our results more rigorously.\\

\subsection{Definitions}~\\

We begin by recalling the definitions of graph edit distance and the edit distance function.

\begin{defn}[Alon-Stav \cite{as2008}; Axenovich-K\'{e}zdy-Martin \cite{akr2008}]
Let $G$ and $H$ be simple graphs on the same labeled vertex set, and let $\Hmc$ be a hereditary property, then
\begin{enumerate}
\item $\dist(G,H)=|E(G)\Delta E(H)|$ is the edit
distance from $G$ to $H$,
\item $\dist(G,\Hmc)=\min\{\dist(G,H):H\in \Hmc\}$ is the edit
distance from $G$ to $\Hmc$ and
\item $\dist(n,\Hmc)=\max\{\dist(G,\Hmc):|G|=n\}$ is the maximum edit distance from the set of all $n$-vertex graphs to the hereditary property $\Hmc$.
\end{enumerate}
\end{defn}

Since $\Hmc$ is by definition closed under isomorphism, vertex labels may be ignored when considering
$\dist(G,\Hmc)$. In fact, $\dist(G,\Hmc)$ could be defined equivalently as the minimum number of edge changes
necessary to make $G$ a member of $\Hmc$.

The limit of the maximum edit distance from an $n$-vertex graph to a hereditary property $\Hmc$ normalized by the
total number of potential edges in an $n$-vertex graph,
$$ d^*_{\Hmc}=\lim_{n\rightarrow\infty}\dist(n,\Hmc)/{\textstyle \binom{n}{2}} , $$
is demonstrated in \cite{as2008} to exist and to be realized asymptotically with high probability by the random
graph $G(n,p^*)$, where $p^*\in [0,1]$ is a probability that depends on $\Hmc$ and is not necessarily unique.

\begin{defn}[\cite{BaloghMartin08}]\label{def:distfunc}
The \textdef{edit distance function} of a hereditary property $\Hmc$ is defined as follows:
$$ \ed_{\Hmc}(p)=\lim_{n\rightarrow\infty}\max\left\{\dist(G,\Hmc) : |V(G)|=n, |E(G)|=\lfloor p{\textstyle \binom{n}{2}}\rfloor\right\}/{\textstyle \binom{n}{2}} . $$
\end{defn}
This function has also been denoted as $g_{\Hmc}(p)$ in, for example,~\cite{BaloghMartin08}.  It was also proven
in~ \cite{BaloghMartin08} that, if $\E$ denotes the expectation, then $$
\ed_{\Hmc}(p)=\lim_{n\rightarrow\infty}\E[\dist(G(n,p),\Hmc)]/{\textstyle \binom{n}{2}} . $$ The limits above
were proven to exist in \cite{BaloghMartin08}, and furthermore, $\ed_{\Hmc}(p)$ is both continuous and
concave down.  As a result, the edit distance function attains a maximum value that is equal to $
d^*_{\Hmc}$. The point, or interval, at which $d^*_{\Hmc}$ is attained is denoted $p^*_{\Hmc}$, and when it is
evident from context, the subscript $\Hmc$ may be omitted from both.

\textbf{Colored regularity graphs (CRGs)} are the building blocks for determining the edit distance function for specific
hereditary properties. We leave the formal definition of CRGs as well as basic facts about them to Section~\ref{sec:CRGs}
for the reader who is unfamiliar with these objects.  A few of the new results given in Section~\ref{sec:results} are
structural in nature and so an understanding of CRGs is useful in order to put the results in full context.\\

\subsection{New Results}~\\
\label{sec:results}

In this paper, we prove the following results for the hereditary properties $\forb(K_{2,3})$ and
$\forb(K_{2,4})$. The case of $K_{2,2}$ is mentioned in Section 5.3 of \cite{MT}.

\begin{thm}\label{mt1}
Let $\Hmc=\forb(K_{2,3})$.  Then $\ed_{\Hmc}(p)=\min\{p(1-p),\frac{1-p}{2}\}$ with $p_\Hmc^*=\frac{1}{2}$ and
$d_\Hmc^*=\frac{1}{4}$.
\end{thm}

\begin{thm}\label{mt2}
Let $\Hmc=\forb(K_{2,4})$.  Then $\ed_{\Hmc}(p)=\min\{p(1-p),\frac{7p+1}{15},\frac{1-p}{3}\}$ with
$p_\Hmc^*=\frac{1}{3}$ and $d_\Hmc^*=\frac{2}{9}$.
\end{thm}

\begin{figure}[htbp]\ \hfill%
\begin{minipage}[t]{2.5in}
\includegraphics[width=2.5in]{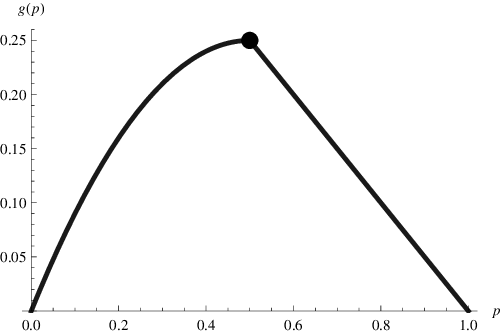}
\caption{Plot of $\ed_{\forb(K_{2,3})}(p)=\min\{p(1-p),(1-p)/2\}$.  The point $(p^*,d^*)=(1/2,1/4)$ is
indicated.}\label{fig:plotk23}
\end{minipage}\ \hfill\ %
\begin{minipage}[t]{2.5in}
\includegraphics[width=2.5in]{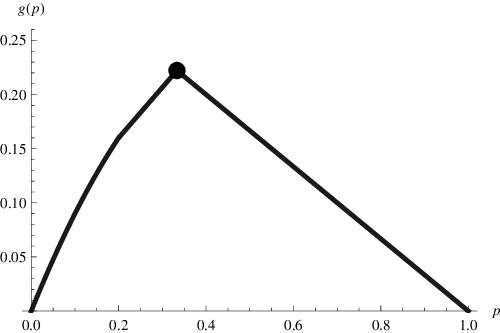}
\caption{Plot of $\ed_{\forb(K_{2,4})}(p)=\min\{p(1-p),(1+7p)/15,(1-p)/3\}$.  The point $(p^*,d^*)=(1/3,2/9)$ is
indicated.} \label{fig:plotk24}
\end{minipage}\ \hfill\ %
\end{figure}

It should be noted that $p_\Hmc^*$ and $d_\Hmc^*$ from Theorem \ref{mt1} could be found using alternative methods
from previous literature as well.  In fact, they are a direct result of Lemma 5.14 in \cite{MT}, as is the value
of $\ed_{\Hmc}(p)$ for $p\geq1/2$ in both theorems. The $p_\Hmc^*$ and $d_\Hmc^*$ values in Theorem \ref{mt2},
however, are not so easily found. The techniques used to prove both theorems also have the potential to yield
some results for $\ed_\Hmc(p)$ when $\Hmc=\text{Forb}(K_{2,t})$ and $t\geq5$, as discussed in Sections \ref{6}
and \ref{7}.

In \cite{MT}, it is established that for $p\geq1/2$ and $\Hmc=\forb(K_{2,t})$, the edit distance function
$\ed_\Hmc(p)=(1-p)/(t-1)$. We extend this result to hold true for $p\geq 2/(t+1)$.
\begin{thm}\label{thm:bdext}
Let $t\geq4$, $p\geq 2/(t+1)$ and $\Hmc=\forb(K_{2,t})$, then $\ed_\Hmc(p)=(1-p)/(t-1)$.
\end{thm}

This extension along with a new CRG construction results in the determination of $d^*_\Hmc$, the maximum value of
the edit distance function, for all odd $t$.  Using the general lower bound in Theorem \ref{thm:genLB} below, we
also demonstrate, via Theorem \ref{thm:oddt}, that this maximum value occurs on a nondegenerate \emph{interval}
of values for $p$.  That is, $p^*_\Hmc$ is not a single value for all odd $t\geq5$.

\begin{thm}\label{thm:genLB}
Let $t\geq 3$ and $p<1/2$.  If $K$ is a black-vertex, $p$-core CRG with white and gray edges such that the gray
edges have neither a $K_{2,t}$ nor a $B_{t-2}$ (as defined in Lemma \ref{lem:forbsub}), then
\begin{equation}\label{eq:genLB} g_K(p)\geq p-\frac{t-1}{4t-5}\left[3p-2+2\sqrt{1-3p+(t+1)p^2}\right] . \end{equation}
\end{thm}

\begin{thm}\label{thm:oddt}
For odd $t\geq5$ and $\Hmc=\forb(K_{2,t})$, $$d_\Hmc^*=1/(t+1)\qquad \mbox{and}\qquad
p^*_\Hmc\supseteq\left[\frac{2t-1}{t(t+1)},\frac{2}{t+1}\right].$$
\end{thm}

For small $p$ and $t$ large enough, we demonstrate a similar result for Zarankiewicz type constructions by
F\"{u}redi \cite{FurediZar} to that discovered in \cite{MT} for $\forb(K_{3,3})$.  Observe that this construction is not used in the ${\rm ed}_{\forb(K_{2,t})(p)}$ bound for $t=3,4$ but is used for $t\geq 5$:
\begin{thm}\label{thm:tg9}
For $\Hmc=\forb(K_{2,t})$, the edit distance function $\ed_\Hmc(p)\leq\frac{t-1+p(2q^2-q(t-1)-2t)}{2(q^2-1)}$ for
any prime power $q$ such that $t-1$ divides $q-1$. \end{thm}

\begin{cor}\label{cor:tg9}
For $t\geq9$, there exists a value $q_0$, so that if $q>q_0$, then
$\frac{t-1+p(2q^2-q(t-1)-2t)}{2(q^2-1)}<p(1-p)$ for some values of $p$, which approach $0$ as $q$ increases. That
is, arbitrarily close to $p=0$, there is some value for $p$ such that $\ed_\Hmc(p)< p(1-p)$.
\end{cor}

A strongly regular graph construction provides the upper bound $\frac{7p+1}{15}$ for $\ed_{\forb(K_{2,4})}(p)$.
Such constructions continue to be relevant for larger $t$ values.

\begin{thm}\label{thm:srgform}
For any $(k,d,\lambda,\mu)$-strongly regular graph, there exists a corresponding CRG, $K$, such that
$$f_K(p)= \frac{1}{k}+\left(\frac{k-d-2}{k}\right)p . $$
\noindent If $\lambda\leq t-3$ and $\mu\leq t-1$, then $K$ forbids $K_{2,t}$ embedding, and when equality holds
for both $\lambda$ and $\mu$,
\begin{equation} \label{eq:SRGfbound} f_K(p)=\frac{t-1}{t-1+d(d+1)}+\left(1-\frac{(d+2)(t-1)}{t-1+d(d+1)}\right)p . \end{equation}
\end{thm}

There is a very close connection between the result in Theorem \ref{thm:genLB} and strongly regular graphs.  If
we take the expression in (\ref{eq:SRGfbound}) and minimize it with respect to $d$ (see expression
(\ref{eqn:theor})), then we obtain the expression on the right-hand side of (\ref{eq:genLB}).  In particular, we
show that if a $(k,d,t-3,t-1)$-strongly regular graph exists, then the corresponding CRG, $K$, has $f_K(p)$
tangent to the curve $p-\frac{t-1}{4t-5}[3p-2+2\sqrt{1-3p+(t+1)p^2}]$ for some value of $p$.  Thus, we obtain the
exact value of $\ed_{\forb(K_{2,t})}$ for that value of $p$.

The following general upper bound arises from a CRG construction involving the second power of cycles.  It is superseded by the strongly regular graph constructions for small values of $t$, but not necessarily for large $t$.
\begin{thm}\label{thm:cycc}
For $\Hmc=\forb(K_{2,t})$, $$\ed_\Hmc(p)\leq\frac{3p+1}{5+t}.$$
\end{thm}
Observe that this bound is better than that of the trivial $\min\{p(1-p),(1-p)/(t-1)$ for $p\in
\left(\frac{t+2-\sqrt{t^2-16}}{2t+10},\frac{3}{2t+1}\right)$ when $t\geq 5$.

Below are known upper bounds for $5\leq t\leq 8$. It should be noted that as our knowledge of existing strongly
regular graphs increases, new upper bounds are also likely to be discovered.

\begin{thm}\label{thm:sum}
Let $\Hmc=\forb(K_{2,t})$.
\begin{itemize}
\item If $t=5$, then
$$\ed_\Hmc(p)\leq\min\left\{p(1-p),\frac{1+75p}{96},\frac{1+26p}{40},\frac{1+5p}{13},\frac{1}{6},\frac{1-p}{4}\right\}.$$
\item If $t=6$, then $$\ed_\Hmc(p)\leq\min\left\{p(1-p),\frac{1+63p}{85},\frac{1+14p}{26},\frac{1+7p}{17},\frac{1+2p}{10},\frac{1-p}{5}\right\}.$$
\item If $t=7$, then
$$\ed_\Hmc(p)\leq\min\left\{p(1-p),\frac{1+124p}{156},\frac{1+76p}{100},\frac{1+44p}{64},\frac{1+31p}{49},\frac{1+20p}{36},\frac{1+5p}{16},\frac{1}{8},\frac{1-p}{6}\right\}.$$
\item If $t=8$, then
$$\ed_\Hmc(p)\leq\min\left\{p(1-p),\frac{1+124p}{156},\frac{1+95p}{125},\frac{1+53p}{76},\frac{1+20p}{36},\frac{1+11p}{25},\frac{1+5p}{16},\frac{3p+1}{13},\frac{1-p}{7}\right\}.$$
\end{itemize}
\end{thm}
We compare these upper bounds to the lower bound in Theorem \ref{thm:genLB} via the figures in Appendix \ref{sec:figs}.~\\

\subsection{CRGs and Background}~\\
\label{sec:CRGs}

To help describe how  $\ed_\Hmc(p)$ may be calculated, some definitions from \cite{as2008} are required.

\begin{defn}[Alon-Stav \cite{as2008}]
A \textdef{colored regularity graph (CRG)}, $K$, is a complete graph with vertices colored black or white, and
with edges colored black, white or gray.
\end{defn}

At times, it may be convenient to refer to the graph induced by edges of a particular color in a CRG, $K$.  We
shall refer to these graphs as the black, white and gray subgraphs of $K$. The investigation of edit distance in
\cite{MT} and in \cite{ThomSurv} uses a different paradigm with an analogous structure called a \emph{type}.
Essentially, our black, white and gray are their red, blue and green, respectively.

\begin{defn}[Alon-Stav \cite{as2008}]
A \textdef{colored homomorphism} from a (simple) graph $H$ to a colored regularity graph $K$ is a mapping
$\phi:V(H)\mapsto V(K)$, which satisfies the following:
\begin{enumerate}
\item  If $uv\in E(H)$, then either $\phi(u)=\phi(v)=m$ and $m$ is colored black, or $\phi(u)\neq\phi(v)$ and the edge $\phi(u)\phi(v)$ is colored black or gray.
\item  If $uv\notin E(H)$, then either $\phi(u)=\phi(v)=m$ and $m$ is colored white, or $\phi(u)\neq\phi(v)$ and the edge $\phi(u)\phi(v)$ is colored white or gray.
\end{enumerate}
\end{defn}

Basically, a colored homomorphism is a map from a  simple graph to a CRG so that black is only associated with
adjacency, white is only associated with nonadjacency and gray is associated with adjacency, nonadjacency or
both. We will refer to a colored homomorphism from a simple graph $H$ to a CRG $K$ as \emph{an embedding of $H$
in $K$}, and we denote the set of all CRGs that only allow the embedding of simple graphs in a hereditary
property $\Hmc$ as $\Kmc(\Hmc)$ or merely $\Kmc$ when $\Hmc$ is clear from the context. Since any hereditary
property may be described by a set of forbidden induced subgraphs, an equivalent description of $\Kmc(\Hmc)$ is
the set of all CRGs that do not permit the embedding of any of the forbidden induced subgraphs associated with
$\Hmc$.  For instance, $\Kmc(\forb(K_{2,3}))$ is the set of all CRGs that do not admit $K_{2,3}$ embedding.

In order to calculate $\ed_\Hmc(p)$, colored regularity graphs are used in \cite{BaloghMartin08} in order to define
the following functions:

\begin{eqnarray}
f_K(p) & = & \frac{1}{k^2}\left[p(|\vw(K)|+2|\ew(K)|)+(1-p)(|\vb(K)|+2|\eb(K)|)\right] \label{eqn:f} \\
g_K(p) & = & \min\{\mathbf{u}^T M_K(p) \mathbf{u}:\mathbf{u}^T\mathbf{1}=1\text{ and
}\mathbf{u}\geq0\}\label{eqn:g} .
\end{eqnarray}
\

 Here $K$ is a CRG with $k$ vertices.  $\vw(K),\ \vb(K),\ \ew(K)$ and $\eb(K)$ represent the sets of white
vertices, black vertices, white edges and black edges in $K$ respectively.  $M_K$ is essentially a weighted
adjacency matrix for $K$ with black vertices and edges receiving weight $1-p$, white vertices and edges receiving
weight $p$ and gray edges receiving weight $0$.  From \cite{as2008}, it is known that
$\ed_\Hmc(p)=\inf_{K\in\Kmc}\{f_K(p)\}=\inf_{K\in\Kmc}\{g_K(p)\}$.  Moreover, Alon and Stav \cite{as2008} show
that if $\chi_B$ is the binary chromatic number of $\Hmc$, then $\ed_\Hmc(1/2)=1/(\chi_B-1)$. Marchant and
Thomason, demonstrate in \cite{MT} that $\ed_\Hmc(p)=\min_{K\in\Kmc}\{g_K(p)\}$. That is, given $p$ there exists
at least one CRG, $K\in\Kmc$, such that $\ed_\Hmc(p)=g_K(p)$.

If we say a CRG, $K$, is a \emph{sub-CRG} of another CRG, $K'$, when $\vw(K)\subseteq\vw(K')$,
$\vb(K)\subseteq\vb(K')$, $\ew(K)\subseteq\ew(K')$ and $\eb(K)\subseteq\eb(K')$, then it may be observed that
$g_K(p)\geq g_{K'}(p)$. Furthermore, as was noted in \cite{MT}, if $g_K(p)= g_{K'}(p)$, then there is no need to
consider both $K$ and $K'$ when attempting to determine $\min_{K\in\Kmc}\{g_K(p)\}$.  Thus, a special subset of
CRGs is defined as follows.

\begin{defn}[Marchant-Thomason \cite{MT}]
A \textdef{$p$-core CRG} is a CRG $K'$ such that for no nontrivial sub-CRG $K$ of $K'$ is it the case that
$g_K(p)=g_{K'}(p)$. In other words, if $K'$ is a $p$-core CRG, and $K$ is a nontrivial sub-CRG of $K'$, then
$g_K(p)>g_{K'}(p)$.
\end{defn}

It can be shown (see \cite{MT}) that a CRG, $K$, is $p$-core if and only if $g_K(p)=\mathbf{x}^TM_K(p)\mathbf{x}$
for a unique vector $\mathbf{x}$ with positive entries summing to $1$.  Any CRG, $K$, that is not $p$-core
contains at least one $p$-core sub-CRG $K'$ so that $g_{K'}(p)=g_{K}(p)$. Thus we could also say that
$$ \displaystyle \ed_\Hmc(p)=\min\{g_K(p):K\in\Kmc\ \text{and}\ K\text{ is }p\text{-core}\} . $$

That is, when looking for CRGs to determine $\ed_\Hmc(p)$, the search may be limited to the important subset of
CRGs, $p$-cores.  This observation is especially helpful for determining lower bounds for the edit distance
function.

To prove the upper bounds for the edit distance function in this paper, we show that for each $p$ there exists a
CRG, $K\in\Kmc$, so that $f_K(p)$ gives the bound for the value of $\ed_\Hmc(p)$. For the lower bounds, we employ
so-called symmetrization, due to Sidorenko~\cite{Sid93}, previously used for the computation of edit distance
functions in \cite{rm2009} and elsewhere such as Marchant and Thomason~\cite{MT}.  We also exploit features of
the graphs $K_{2,t}$ and the concavity of the edit distance function to demonstrate that for no $p$-core CRG,
$K\in\Kmc$, can $g_K(p)$ be less than the value in the theorem.

By the continuity of the edit distance function, if we know the value of the function on an open interval, then
we also know the value on its closure.  Hence, for convenience, most of our proofs will only address the value of
the function on the interior of a given interval.  We also note that, in \cite{rm2009}, whenever the edit
distance function of a hereditary property is computed, there is an attempt to determine all of the $p$-core CRGs
that achieve the value of the edit distance function.  In this paper we only concern ourselves with the value of
the edit distance function itself and do not address the issue of multiple defining constructions.\\

\subsection{The Zarankiewicz problem and strongly regular graphs}\label{1.3}~\\

One reason for our interest in the edit distance function for Forb$(K_{2,t})$ is its relation to the Zarankiewicz
problem.  This problem addresses the question of how many edges a graph can have before it must
contain a $K_{s,t}$ subgraph for fixed $s$ and $t$. In an intriguing result from \cite{MT}, a
construction from Brown \cite{Brown} for $K_{3,3}$-free graphs is applied to construct an infinite set of new
CRGs that improve upon the previously known bounds for $\ed_{\text{Forb}(K_{3,3})}(p)$ on certain intervals for
arbitrarily small $p$.

Marchant and Thomason \cite{MT} establish that it is sufficient to consider only $p$-core CRGs for which the gray
subgraph has neither a $K_3$ nor a $K_{3,3}$.  Brown's constructions are not $K_3$-free, but a bipartite graph
can be created from the construction that has no copy of $K_{3,3}$.  Similarly, for $K_{2,t}$, it is sufficient
to have no gray $K_{2,t}$ or book $B_{t-2}$ to forbid $K_{2,t}$ embedding, where the graph $B_{t-2}$ is a
``book'' as defined in \cite{CFG} and is defined precisely in Lemma \ref{lem:forbsub}.

Although the Brown constructions show that the edit distance function for $\forb(K_{3,3})$ is strictly less than
$p(1-p)/(1+p)$ for sufficiently small $p$, known constructions for dense $K_{2,t}$-free graphs do not play a role
in the computation of the edit distance function for $\forb(K_{2,3})$ or $\forb(K_{2,4})$ in the same way.
However, the edit distance function for $\forb(K_{2,4})$ is achieved over the interval $[1/5,1/3]$ by a
construction formed from a strongly regular graph, namely a generalized quadrangle, often denoted ${\rm
GQ}(2,2)$, and similar results to those from the Brown constructions do reemerge when $t\geq 9$ and F\"{u}redi's
constructions from \cite{FurediZar} are considered.\\

\subsection{Organization}~\\

In Section \ref{2}, we discuss some results from \cite{MT} and \cite{rm2009} for the edit distance function and
how they may be applied to the problem of determining the function for Forb$(K_{2,t})$. We then proceed to some
general results and observations in Section \ref{3} that will be useful throughout the paper. Sections \ref{4}
and \ref{5} contain the proofs of our results for Forb$(K_{2,3})$ and Forb$(K_{2,4})$, respectively. Section
\ref{6} addresses the proofs of Theorems \ref{thm:bdext}, \ref{thm:genLB} and \ref{thm:oddt}.  In Section
\ref{7}, we present several new CRG constructions that yield upper bounds for $\ed_{\forb(K_{2,t})}(p)$ in
general.  The remaining sections are reserved for conclusions and acknowledgements.

\section{Applications of past results to Forb$(K_{2,t})$}\label{2}

If a CRG is $p$-core, one can say some interesting things about its overall structure. From \cite{MT}, for
instance, we have the following useful result.
\begin{thm}[Marchant-Thomason \cite{MT}]\label{edges}
If $K$ is a $p$-core CRG, then all edges of $K$ are gray except
\begin{itemize}
\item if $p<1/2$, some edges joining two black vertices might be white or
\item if $p>1/2$, some edges joining two white vertices might be black.
\end{itemize}
\end{thm}

The CRGs with all edges gray are useful in bounding the edit distance function, as we see in \cite{rm2009}.
\begin{defn}  Let $K(w,b)$ denote the CRG with $w$ white vertices, $b$ black vertices and only gray edges.  In particular:
\begin{enumerate}
\item Let $K(1,1)$ be the CRG consisting of a white and black vertex joined by a gray edge.
\item Let $K(0,t-1)$ be the CRG consisting of $t-1$ black vertices all joined by gray edges.
\end{enumerate}
\end{defn}

Theorem \ref{thm:k22} settled the case of $K_{2,2}$, and Theorem \ref{thm:trivbds} permits us to focus on
black-vertex CRGs.

\begin{thm}[Marchant-Thomason \cite{MT}]\label{thm:k22}
Let $\Hmc=$Forb$(K_{2,2})$.  Then $\ed_{\Hmc}(p)=g_{K(1,1)}(p)=p(1-p)$ with $p_\Hmc^*=\frac{1}{2}$ and
$d_\Hmc^*=\frac{1}{4}$.
\end{thm}

\begin{thm}[Marchant-Thomason \cite{MT}]\label{thm:trivbds}
Let $\Hmc=\forb(K_{2,t})$, $t>2$.  Then
\begin{enumerate}
\item For $p>\frac{1}{2}$, $\ed_{\Hmc}(p)=g_{K(0,t-1)}(p)=\frac{1-p}{t-1}$ and
\item For $p\leq\frac{1}{2}$, either
\begin{itemize}
\item $\ed_\Hmc(p)=\min\{g_{K(1,1)}(p),g_{K(0,t-1)}(p)\},$ or
\item $\ed_\Hmc(p)=g_K(p)<\min\{g_{K(1,1)}(p),g_{K(0,t-1)}(p)\},$ where $K$ is a $p$-core CRG with only black vertices and, consequently, no black edges.
\end{itemize}
\end{enumerate}
\end{thm}

The following lemma is about the structure of the $p$-core CRGs described in the second part of Theorem
\ref{thm:trivbds}.  It was originally observed in Example 5.16 of \cite{MT}.  The proof is a straightforward case
analysis.

\begin{lem}[Marchant-Thomason \cite{MT}]\label{lem:forbsub} A CRG, $K$, on all black vertices with only white and gray edges forbids $K_{2,t}$
embedding if and only if its gray subgraph contains no $K_{2,t}$ or $\T$ as a subgraph, where $\T$ is a book as
described in \cite{CFG}. That is, the graph $\T$ is defined to be the graph consisting of $t-2$ triangles that
all share a single common edge.
\end{lem}

As demonstrated in \cite{MT}, for a $p$-core CRG, $K$, there is a \emph{unique} vector $\xb$ so that
$g_K(p)=\xb^T M_K(p) \xb$.
\begin{defn}[Marchant-Thomason \cite{MT}]
For a $p$-core CRG $K$ with optimal weight vector $\xb$, the entry of $\xb$ corresponding to a vertex, $v\in
V(K)$, is denoted by $\xb(v)$.  This is the \textdef{weight of $v$}, and the function $\xb(v)$ is the
\textdef{optimal weight function}.
\end{defn}

With this in mind, we have two propositions from \cite{rm2009}, which follow easily from \cite{MT}.
\begin{prop}[\cite{rm2009}]\label{prop:dgformula}
Let $K$ be a p-core CRG with all vertices black.  Then for any $v\in V(K)$ and optimal weighting $\xb$,
$d_G(v)=\frac{p-g_K(p)}{p}+\frac{1-2p}{p}\xb(v)$, where $d_G(v)$ is the sum of the weights of the vertices
adjacent to $v$ via a gray edge.
\end{prop}

\begin{prop}[\cite{rm2009}]\label{prop:xbound}  Let $K$ be a p-core CRG with all
vertices black, then for $p\in [0,1/2]$ and optimal weighting $\xb$,
$$\xb(v)\leq\frac{g_K(p)}{1-p},\ \forall v\in V(K).$$
\end{prop}

Because of Theorem \ref{thm:trivbds}, we may restrict our attention to those CRGs, $K$, for which $g_K(p)\leq
p(1-p)$.  As a result, Proposition \ref{prop:dgformula} gives the lower bound $d_G(v)\geq
p+\frac{1-2p}{p}\xb(v)$. Meanwhile, Proposition \ref{prop:xbound} restricts the optimal weights of all vertices
in $K$ to be no more than $p$. These two restrictions are useful when attempting to prove lower bounds for
$\ed_{\text{Forb}(K_{2,t})}(p)$.

\section{Preliminary results and observations}\label{3}

We begin with some notation used throughout the paper.
\begin{defn}
Let $K$ be a black-vertex, $p$-core CRG with $g_K(p)\leq p(1-p)$ and optimal weight function $\xb$:
\begin{itemize}
\item $N_G(v)=\{y\in V(K):\ vy\in EG(K)\}$,
\item $u_0$ is a fixed vertex in $K$ such that $\xb(u_0)\geq\xb(v)$,
 for all $v\in V(K)$, and $x=\xb(u_0)$ is its weight,
\item  $U=N_G(u_0)$ and $|U|=\ell$,
\item $u_1$ is a fixed vertex with maximum weight in $U$, and $x_1=\xb(u_1)$,
\item $W$ is the set of all vertices in $K$ that are neither $u_0$, nor contained in $U$; or equivalently, $W$ is the set of all vertices in the white neighborhood of $u_0$, and
\item $\xb(S)=\sum_{y:y\in S}\xb(y)$ for some set $S\subseteq V(K)$.
\end{itemize}
\end{defn}

\begin{figure}[htbp]
\begin{center}
\includegraphics[width=60mm]{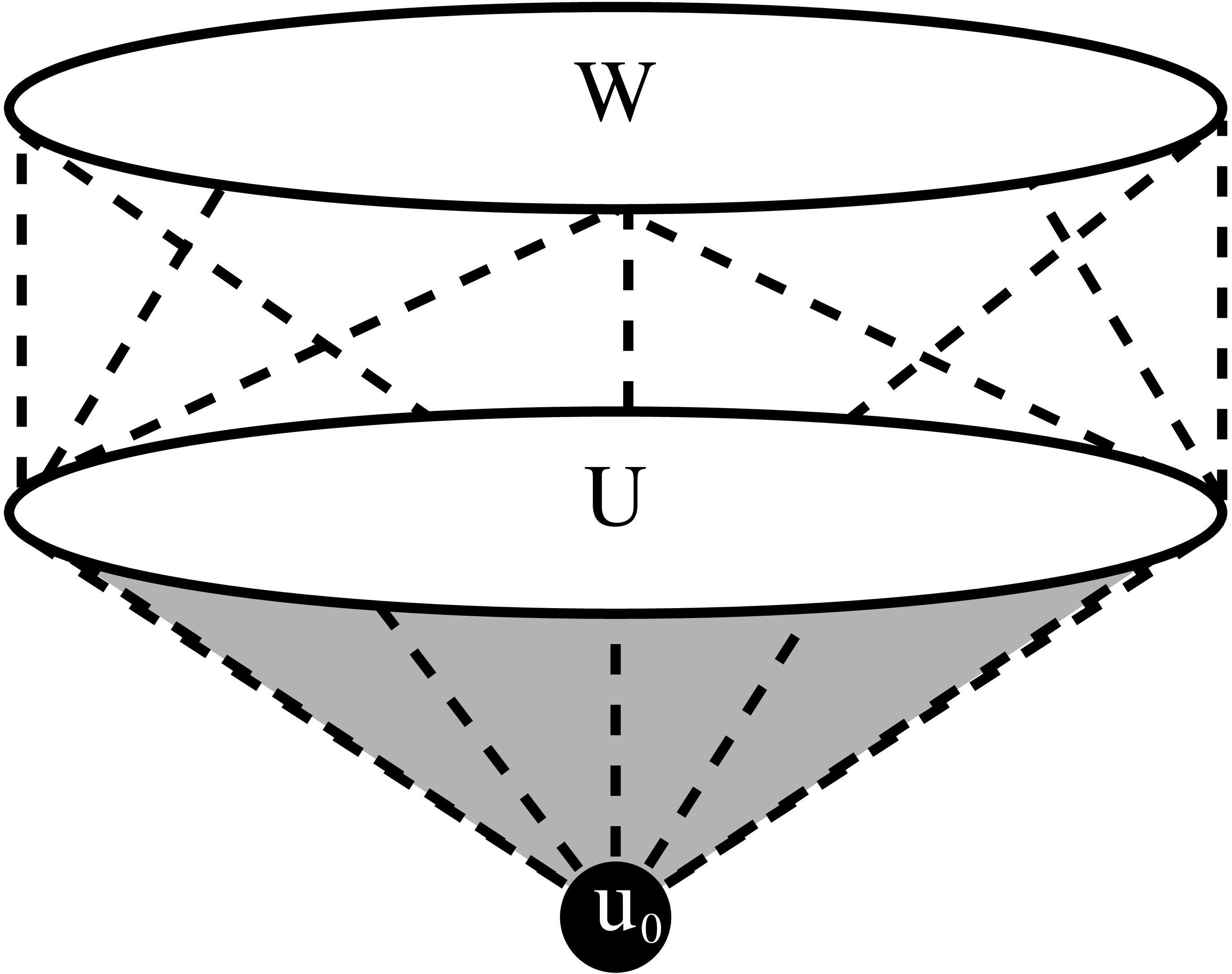}

\caption{A partition of the vertices in a black-vertex, $p$-core CRG, $K$. Dashed lines and gray background
represent gray edges. White edges are omitted, as are edges within subsets.}\label{partition}
\end{center}
\end{figure}

Partitioning the vertices in a black-vertex, $p$-core CRG that forbids a $K_{2,t}$ embedding into the three sets
$\{u_0\},\ U$ and $W$ as seen in Figure \ref{partition}, illustrates some interesting features of its optimal
weight function when the gray neighborhoods of these vertices are examined. One such feature is the upper bounds
in Proposition \ref{starting prop} for $x_1$.

\begin{prop}\label{starting prop}
Let $K\in\left[\Kmc(\forb(K_{2,3}))\cup \Kmc(\forb(K_{2,4}))\right]$ be a black-vertex, $p$-core CRG. If either
$p<1/3$ or both $p<1/2$ and the gray sub-CRG of $K$ is triangle-free, then
$$ x_1\leq x\qquad \mbox{and}\qquad x_1\leq p-x $$
\noindent where $x=\xb(u_0)$ is the maximum weight of a vertex in $K$, and $x_1=\xb(u_1)$ is the maximum weight
of a vertex in that vertex's gray neighborhood.
\end{prop}

\begin{proof}
The inequality $x_1\leq x$ follows directly from definitions of $x_1$ and $x$, since $x$ is the greatest weight
in $K$. To justify the inequality $x_1\leq p-x$, we break the problem into two cases:\newline

\noindent \textbf{Case 1:} $u_0$ and $u_1$ have no common gray neighbor.\\

Recall that $u_1$ is a vertex with maximum weight in the gray neighborhood of $u_0$, a vertex with maximum weight
in all of $K$, and assume that $x+x_1> p$. Then applying Proposition \ref{prop:dgformula} and Theorem
\ref{thm:trivbds},

$$ d_G(u_0)+d_G(u_1) \geq \left[p+\frac{1-2p}{p}x\right]+\left[p+\frac{1-2p}{p}x_1\right] = 2p+\left(\frac{1-2p}{p}\right)(x+x_1) > 2p+(1-2p) . $$

This is a contradiction because in Case 1, $N_G(u_0)\cap N_G(u_1)=\emptyset$.  Thus, $d_G(u_0)+d_G(u_1)\leq1$, since the sum of the weights of the vertices in $K$
must be $1$.

This completes the proof of Proposition~\ref{starting prop} for $K\in\Kmc(\forb(K_{2,3}))$, since, in this case,
no $K\in\Kmc$ contains a gray triangle.  So we may assume that $K\in\Kmc(\forb(K_{2,4}))$.\newline

\noindent \textbf{Case 2:} $u_0$ and $u_1$ have a common gray neighbor and $p<1/3$. \\

In this case, $u_1$ has a single neighbor $u_2$ in $U$ because any more such neighbors would result in a gray
book $B_2$.  Furthermore, we note that in order to avoid a gray book $B_2$, the common neighborhood of $u_1$ and
$u_2$ in $W$ must be empty.  Consequently, $d_G(u_1)+d_G(u_2)\leq \xb(W)+2x+x_1+\xb(u_2)$.

Applying similar reasoning to that in Case 1,
$$ d_G(u_0)+d_G(u_1)+d_G(u_2) \geq \left[p+\frac{1-2p}{p}x\right]+\left[p+\frac{1-2p}{p}x_1\right]+
\left[p+\frac{1-2p}{p}\xb(u_2)\right] . $$

So,
\begin{eqnarray*}
d_G(u_0)+\left(\xb(W)+2x+x_1+\xb(u_2)\right) & \geq &
\left[p+\frac{1-2p}{p}x\right]+\left[p+\frac{1-2p}{p}x_1\right]+
\left[p+\frac{1-2p}{p}\xb(u_2)\right] \\
\xb(U)+\left(\xb(W)+2x+x_1+\xb(u_2)\right) & \geq & 3p+\frac{1-2p}{p}\left(x+x_1+\xb(u_2)\right) \\
\xb(U)+\xb(W)+x & \geq & 3p+\frac{1-3p}{p}\left(x+x_1+\xb(u_2)\right) \\
1 & \geq & 3p+\frac{1-3p}{p}\left(x+x_1+\xb(u_2)\right). \\
\end{eqnarray*}

With $p<1/3$ and $x+x_1\geq p$, we have a contradiction.
\end{proof}

Applying the pigeon-hole principle, we also have the following lower bound for $\ell$:

\begin{fact}\label{ph}
In a CRG, if $u_0$ is a vertex with maximum weight, $x=\xb(u_0)$, the maximum weight in the gray neighborhood of
$u_0$ is $x_1$, and the order of the gray neighborhood of $u_0$ is $\ell$, then $\ell\geq d_G(u_0)/x_1\geq
d_G(u_0)/x$.
\end{fact}

While simple, when combined with Propositions \ref{prop:dgformula} and \ref{starting prop} along with the
observation that $\xb(u_0)+\xb(U)+\xb(W)=1$, this fact forces a delicate balance between the weights of the
vertex $u_0$, the vertices in $U$, and the vertices in $W$.

\section{Proof of Theorem \ref{mt1}}\label{4}
In this section, we establish the value of $\ed_{\forb(K_{2,3})}(p)$ for $p\in(0,1/2)$, determining the entire
function via continuity and Theorem \ref{thm:trivbds}, which gives that $\ed_{\forb(K_{2,3})}(p)=(1-p)/2$ for
$p\in[1/2,1]$.

For the following discussion, we assume that $K$ is a $p$-core CRG on all black vertices into which $K_{2,3}$ may
not be embedded and that $g_K(p)\leq p(1-p)$.  The following lemma yields a useful restriction of the order of
$U$.

\begin{lem}\label{upperbounds for U}  Let $K$ be a black-vertex, $p$-core CRG with $p\in(0,1/2)$, no gray triangles,
no gray $K_{2,3}$ and $g_K(p)\leq p(1-p)$.  If $u_0$ is a vertex of maximum weight, $x$, in $K$, and
$\ell=|N_G(u_0)|$, then
$$ \ell\leq \frac{2(1-x)-\frac{1}{p}d_G(u_0)}{p-x} . $$
\end{lem}

\begin{proof}
Let $u_1,\ldots,u_\ell$ be an enumeration of the vertices in $U$, the gray neighborhood of $u_0$. Observe that
$K$ cannot contain a $K_3$ with all gray edges, and so $U$ contains no gray edges. Therefore, with the
exception of $u_0$, the entire gray neighborhood of each $u_i$ is contained in $W$. Furthermore, if any three
vertices in $U$ had a common gray neighbor in $W$, then $K$ would contain a gray $K_{2,3}$. That is, each vertex
in $W$ is adjacent to at most $2$ vertices in $U$ via a gray edge. Applying these observations,
$$ \sum_{i=1}^\ell(d_G(u_i)-x) \leq 2\xb(W) . $$

Using Proposition \ref{prop:dgformula} and the assumption that $\frac{p-g_K(p)}{p}\geq p$,
$$ \sum_{i=1}^\ell\left(p-x+\frac{1-2p}{p}\xb(u_i)\right)\leq 2\xb(W) . $$
The fact that $\xb(W)=1-x-d_G(u_0)$, gives
\begin{eqnarray*}
\ell(p-x)+\frac{1-2p}{p}d_G(u_0) & \leq & 2\left(1-x-d_G(u_0)\right) \\
\ell(p-x) & \leq & 2-2x-\frac{1}{p}d_G(u_0) \\
\ell&\leq& \frac{2(1-x)-\frac{1}{p}d_G(u_0)}{p-x} .
\end{eqnarray*}
\end{proof}

The following technical lemma is an important tool in the proof of the Theorem.

\begin{lem}\label{k23 thrm}
Let $K$ be a black-vertex, $p$-core CRG for $p\in(0,1/2)$ with no gray triangles, no gray $K_{2,3}$ and
$g_K(p)\leq p(1-p)$.  If $x$ and $x_1$ are defined as in Proposition \ref{starting prop}, then
$$ \left[p+\frac{1-2p}{p} x\right]\left[\frac{1}{x_1}+\frac{1}{p(p-x)}\right] \leq \frac{2(1-x)}{p-x} . $$
\end{lem}

\begin{proof}
By Fact \ref{ph}, $\ell\geq\frac{d_G(u_0)}{x_1}$, and by Lemma \ref{upperbounds for U}, $\ell\leq
\frac{2(1-x)-\frac{1}{p}d_G(u_0)}{p-x}$.  Therefore,
$$\frac{d_G(u_0)}{x_1}\leq\frac{2(1-x)-\frac{1}{p}d_G(u_0)}{p-x}.$$

After combining the $d_G(u_0)$ terms we get,
\begin{eqnarray*}
d_G(u_0)\left[\frac{1}{x_1}+\frac{1}{p(p-x)}\right]&\leq&\frac{2(1-x)}{p-x},
\end{eqnarray*}
and then applying Proposition \ref{prop:dgformula},
\begin{eqnarray*}
\left[p+\frac{1-2p}{p}x\right]\left[\frac{1}{x_1}+\frac{1}{p(p-x)}\right]&\leq&\frac{2(1-x)}{p-x}.
\end{eqnarray*}
\end{proof}

We now turn to the proof of the main theorem for this section.

\begin{proof}[Proof of Theorem \ref{mt1}]
Let $p\in (0,1/2)$, and $K$ be a black-vertex, $p$-core CRG with $g_K(p)<p(1-p)$ and no gray triangle (i.e., the
book $B_1$) or gray $K_{2,3}$.

With the above assumptions, we will show that there is no possible value for $x$, the value of the largest
vertex-weight.  To do so, we break the problem into 2 cases: $x\geq\frac{p}{2}$ and $x<\frac{p}{2}$.\newline

\noindent\textbf{Case 1:} $x\geq p/2$.\\

We start with the inequality from Lemma \ref{k23 thrm},
\begin{eqnarray*}
\left[p+\frac{1-2p}{p}x\right]\left[\frac{1}{x_1}+\frac{1}{p(p-x)}\right]&\leq&\frac{2(1-x)}{p-x},
\end{eqnarray*}
and apply the bound $x_1\leq p-x$ from Proposition \ref{starting prop} to get
\begin{eqnarray*}
\left[p+\frac{1-2p}{p}x\right]\left[\frac{1}{p-x}+\frac{1}{p(p-x)}\right]&\leq&\frac{2(1-x)}{p-x}.\\
\end{eqnarray*}
From Proposition \ref{prop:xbound}, $p-x>0$, and so
\begin{eqnarray*}
\left[p+\frac{1-2p}{p}x\right]\left[1+\frac{1}{p}\right] & \leq & 2(1-x)\\
x\left(\frac{1-p}{p^2}\right) & \leq & 1-p\\
x & \leq & p^2,\\
\end{eqnarray*}
\noindent a contradiction, since $\frac{p}{2}>p^2$ for $p\in(0,1/2)$.\newline

\noindent \noindent \textbf{Case 2:} $x< p/2$.\\

We again apply Lemma \ref{k23 thrm}, only now we employ the trivial bound $x_1\leq x$ from Proposition
\ref{starting prop}:

\begin{eqnarray*}
\left[p+\frac{1-2p}{p}x\right]\left[\frac{1}{x}+\frac{1}{p(p-x)}\right]&\leq&\frac{2(1-x)}{p-x}\\
\left[p+\frac{1-2p}{p}x\right]\left[p(p-x)+x\right]&\leq&2px(1-x)\\
(4p^2-3p+1)x^2-(3p^3)x+p^4&\leq& 0.\\
\end{eqnarray*}

Observe that $4p^2-3p+1$ is always positive, and therefore the parabola $(4p^2-3p+1)x^2-(3p^3)x+p^4$, in the variable $x$, is concave up, so the range of $x$ values for which this inequality is satisfied is $x\in \left[x',x''
\right]$ where\newline

$$x'=\frac{3p^3-\sqrt{-4p^4+12p^5-7p^6}}{2(1-3p+4p^2)}\qquad \mbox{and}\qquad
x''=\frac{3p^3+\sqrt{-4p^4+12p^5-7p^6}}{2(1-3p+4p^2)}.$$\newline

If $p<(6-2\sqrt{2})/7$, then neither $x'$ nor $x''$ is real, and so the inequality is never satisfied. For
$p\in\left[\frac{6-2\sqrt{2}}{7},\frac{1}{2}\right)$, routine calculations show that $\frac{p}{2}<x'$, a
contradiction to the assumption that $x<\frac{p}{2}$.

Hence, there is no possible value for $x$ if $\ed_{\text{Forb}(K_{2,3})}(p)<p(1-p)$, so the proof is complete.
\end{proof}

\section{Proof of Theorem \ref{mt2}}\label{5}
This section addresses the case of $\ed_{\forb(K_{2,4})}(p)$.\\

\subsection{Upper bounds}\label{5.1}~\\

Recall that from Theorem \ref{thm:trivbds} we already know that $\ed_{\forb(K_{2,4})}(p) \leq
\min\{p(1-p),\frac{1-p}{3}\}$. For the remaining upper bound, we turn to strongly regular graphs:
\begin{defn}
A \textdef{$(k,d,\lambda,\mu)$-strongly regular graph} is a graph on $k$ vertices such that each vertex has
degree $d$, each pair of adjacent vertices has exactly $\lambda$ common neighbors, and each pair of nonadjacent
vertices has exactly $\mu$ common neighbors.
\end{defn}

\begin{lem}\label{lem:SRG UB}
Let $\Hmc=\forb(K_{2,t})$.  If there exists a $(k,d,\lambda,\mu)$-strongly regular graph with $\lambda\leq t-3$
and $\mu\leq t-1$, then
$$ \ed_{\Hmc}(p)\leq \frac{1}{k}+\frac{k-d-2}{k}p . $$
\end{lem}

\begin{proof}
Let $G$ be the aforementioned strongly regular graph.  We construct a CRG, $K$, on $k$ black vertices with gray
edges in $K$ corresponding to adjacent vertices in $G$ and white edges in $K$ corresponding to nonadjacent
vertices in $G$.

No pair of adjacent vertices has $t-2>\lambda$ common neighbors, so there is no book $B_{t-2}$ in the gray
subgraph, and no pair of vertices has $t>\mu,\lambda$ common neighbors, so there is no $K_{2,t}$ in the gray
subgraph. Thus, by Lemma \ref{lem:forbsub}, $K_{2,t}\not\mapsto K$.  Furthermore,
$$ f_K(p)=\frac{1}{k^2}\left[(1-p)k+2p\left(\binom{k}{2}-\frac{dk}{2}\right)\right] =\frac{1}{k}+\frac{k-d-2}{k}p . $$
\end{proof}

In fact, there is a $(15,6,1,3)$-strongly regular graph \cite{Brouwer}.  It is a so-called ``generalized
quadrangle,'' ${\rm GQ}(2,2)$.  As a result,
$$ \ed_{\forb(K_{2,4})}(p) \leq \min\left\{p(1-p),\frac{1+7p}{15},\frac{1-p}{3}\right\} . $$

A list of known strongly regular graphs and their parameters has been compiled by Andries Brouwer \cite{Brouwer}.
Their implications for $\ed_{\text{Forb}(K_{2,t})}(p)$ in general are explored further in Section \ref{7.1}.\\

\subsection{Lower bounds}~\\

Because the edit distance function is both continuous and concave down, it is sufficient to verify that
$\ed_{\forb(K_{2,4})}(p)\geq p(1-p)$ for $p\in (0,1/5)$ and that $\ed_{\forb(K_{2,4})}(p)\geq (1-p)/3$ for $p\in
(1/3,1/2)$.  This is because the line determined by the bound $\frac{1+7p}{15}$ passes through the points
$(1/5,4/25)$ and $(1/3,2/9)$. Furthermore, by Theorem \ref{thm:trivbds}, we need only consider CRGs that have
black vertices and white and gray edges.

Lemmas \ref{p large} and \ref{p small} settle the cases where $p\in (1/3,1/2)$ and where $p\in (0,1/5)$,
respectively.

\begin{lem}\label{p large}
Let $p\in (1/3,1/2)$.  If $K$ is a black-vertex, $p$-core CRG that does not contain a gray book $B_2$ or a gray
$K_{2,4}$, then $g_K(p)\geq\frac{1-p}{3}$, with equality occurring only if $K$ is a gray triangle (i.e.,
$K\approx K(0,3)$).
\end{lem}
\begin{proof}

See Appendix~\ref{sec:lemproofs}.

\end{proof}

Below are Lemma \ref{p small} and two propositions that are necessary and used in several cases for the proof of the lemma. 
\begin{prop}\label{p small 1}
Let $p\in (0,1/2)$, and let $K$ be a black-vertex, $p$-core CRG with no gray book $B_2$ and no gray $K_{2,4}$. If
$g=g_K(p)$, $U=N_G(u_0)$, $\ell=|U|$ and $U_1\subseteq U$ is the set of vertices in $U$ that are incident to a
gray edge in $U$, then
$$ \ell\left(\frac{p-g}{p}-x\right)\leq 3-3x-\frac{1+p}{p}\xb(U)+\xb(U_1)\leq3-3x-\frac{1}{p}\xb(U) . $$
\end{prop}

\begin{proof}
See Appendix~\ref{sec:lemproofs}.
\end{proof}~\\

\begin{prop}\label{p small 2}
Let $p\in (0,1/2)$, and let $K$ be a black-vertex, $p$-core CRG with no gray book $B_2$ and no gray $K_{2,4}$. If
$g_K(p)\leq p(1-p)$, then both
$$ p\geq\frac{9-4\sqrt{3}}{11}\qquad\mbox{and}\qquad x\geq\frac{p^2}{2(1-3p+5p^2)}\left[1+3p-\sqrt{-3+18p-11p^2}\right]\geq\frac{1}{25} . $$
\end{prop}

\begin{proof}

See Appendix~\ref{sec:lemproofs}.

\end{proof}

\begin{lem}\label{p small}
Let $p\in (0,1/5)$.  If $K$ is a black-vertex, $p$-core CRG that does not contain a gray book $B_2$ or a gray
$K_{2,4}$, then $g_K(p)> p(1-p)$.
\end{lem}

\begin{proof}

See Appendix~\ref{sec:lemproofs}.

\end{proof}

\section{Proofs of Theorems \ref{thm:bdext}, \ref{thm:genLB} and \ref{thm:oddt}}\label{6}

In this section we extend the generally known interval for $\ed_{\forb(K_{2,t})}(p)$ from $p\in[1/2,1]$ to
$p\in[\frac{2}{t+1},1]$. With a new CRG construction, this extension is sufficient to determine $d^*_\Hmc$ and a
subset of $p^*_\Hmc$ for odd $t$. Subsection \ref{6.1} contains the proof of Theorem \ref{thm:bdext}, while the
remaining subsections address Theorems \ref{thm:genLB} and \ref{thm:oddt}.

\subsection{An extension of the known interval for $K_{2,t}$}\label{6.1}
\begin{proof}[Proof of Theorem \ref{thm:bdext}]
Let $K$ be a $p$-core CRG for $p\in [\frac{2}{t+1},1]$ that does not permit $K_{2,t}$ embedding for
$t\geq 5$. If we assume that $g_K(p)<g_{K(0,t-1)}(p)=(1-p)/(t-1)$, then by Theorem \ref{thm:trivbds}, $K$ has
only black vertices and no black edges.

Again, we partition the vertices of $K$ into three sets $\{u_0\}$, $U=\{u_1,\ldots,u_\ell\}$ and $W$, where $u_0$
is a fixed vertex with maximum weight $x$; $U$ is the set of all vertices in the gray neighborhood of $u_0$ with
$u_1$ a vertex of maximum weight $x_1$ in $U$; and $W$ is the set of all remaining vertices, or those vertices
adjacent to $u_0$ via white edges.  Finally, let $d_G(u_i)$ signify the sum of the weights of all vertices in the
gray neighborhood of $u_i$.

Then by Lemma \ref{lem:forbsub}, the total weight of the vertices in $W$ is at least
$$d_G(u_1)-(t-3)x_1-x,$$

\noindent since no vertex in $U$ can be adjacent to more than $t-3$ other vertices in $U$ without forming a book
$B_{t-2}$ gray subgraph with $u_0$.  Thus,
$$x+d_G(u_0)+\left[d_G(u_1)-(t-3)x_1-x\right]\leq1.$$
Applying Proposition \ref{prop:xbound} and letting $g_K(p)=g$,
\begin{eqnarray*}
2\left(\frac{p-g}{p}\right)+\frac{1-2p}{p}x+\left[\frac{1-2p}{p}-(t-3)\right]x_1 & \leq & 1 \\
2(p-g)-p+(1-2p)x & \leq & \left[p(t-1)-1\right]x_1 \\
2(p-g)-p+(1-2p)x & \leq & \left[p(t-1)-1\right]x \\
p-2g & \leq & \left[p(t+1)-2\right]x.
\end{eqnarray*}
Since $p\geq\frac{2}{t+1}$ and $x\leq\frac{g}{1-p}$ by Proposition \ref{prop:xbound},
\begin{eqnarray*}
p-2g & \leq & \left[p(t+1)-2\right]\frac{g}{1-p} \\
\frac{1-p}{t-1} & \leq & g.
\end{eqnarray*}

By Theorem \ref{thm:trivbds}, $g\leq\frac{1-p}{t-1}$, for $p\in[\frac{2}{t+1},1]$, so
$\ed_{\forb(K_{2,t})}(p)=\frac{1-p}{t-1}$.
\end{proof}

We will now show that this result is enough to determine the maximum value of $\ed_{\forb(K_{2,t})}(p)$ for odd
$t$.

\subsection{A construction for odd $t$}
\begin{prop}\label{prop:t+1}
Let $\Hmc=\forb(K_{2,t})$ for odd $t$.  Then $\ed_\Hmc(p)\leq1/(t+1)$.
\end{prop}

\begin{proof}
Let $K$ be the CRG consisting of $t+1$ black vertices with white subgraph forming a perfect matching and all
other edges gray.  The CRG, $K$, does not contain a gray $K_{2,t}$ or book $B_{t-2}$, and so by Lemma
\ref{lem:forbsub}, $K$ forbids a $K_{2,t}$ embedding.

The CRG, $K$, contains exactly $(t+1)/2$ white edges, so by Equation (\ref{eqn:f}),
$$f_K(p)=\frac{1}{(t+1)^2}\left[p\left(2\cdot\frac{t+1}{2}\right)+(1-p)(t+1)\right]=\frac{1}{t+1}.$$

\noindent Therefore, $\ed_\Hmc(p)\leq 1/(t+1)$.
\end{proof}
Since by Theorem \ref{thm:bdext}, $\ed_{\forb(K_{2,t})}(\frac{2}{t+1})=\frac{1}{t+1}$, and by Proposition
\ref{prop:t+1}, $\ed_{\forb(K_{2,t})}\leq \frac{1}{t+1}$, we have that
$d^*_{\forb(K_{2,t})}=\frac{1}{t+1}$ for odd $t\geq5$.~\\

\subsection{A general lower bound for $t$}~\\

We conclude this section by determining a general lower bound for the edit distance function of $\forb(K_{2,t})$.
It is the lower bound from Theorem \ref{eq:genLB}, and it allows us to make the claim in Theorem \ref{thm:oddt}
that, in the case of odd $t$, there is a nondegenerate interval $p^*_\Hmc$ that achieves the maximum value of the
function.

\begin{proof}[Proof of Theorem \ref{thm:genLB}]
Here we use the standard bounds from Propositions \ref{prop:dgformula} and \ref{prop:xbound}.  Let $g=g_K(p)$,
where $K$ is a black-vertex, $p$-core CRG, and let $N_G(v)$ denote the gray neighborhood of a given vertex $v$ in
$K$. Then if $u_1,\ldots,u_\ell$ are the vertices in the gray neighborhood, $U$, of a fixed vertex of maximum
weight, $u_0$,
$$\sum_{i=1}^{\ell}\left[d_G(u_i)-x-\x\left(N_G(u_i)\cap N_G(u_0)\right)\right] \leq(t-1)(1-x-d_G(u_0)).$$

The left-hand side of this inequality calculates the weight of the total gray neighborhood of each vertex in $U$
that must be contained in $W$, the set of all vertices not in $U$ or $u_0$.  On the right-hand side we make use
of the facts that $\xb(W)=1-x-d_G(u_0)$ and that no vertex in $W$ may be adjacent to more than $(t-1)$ vertices
in $U$ without violating Lemma \ref{lem:forbsub} by forming a gray $K_{2,t}$ with $u_0$.  Thus, applying
Proposition \ref{prop:dgformula},

$$ \sum_{i=1}^{\ell}\left[\frac{p-g}{p}-x+\frac{1-2p}{p}\x(u_i)\right] -\sum_{i=1}^{\ell}\x\left(N_G(u_i)\cap N_G(u_0)\right)\leq(t-1)(1-x-d_G(u_0)).$$

Again considering Lemma \ref{lem:forbsub} reveals that no vertex $u_i\in U$ can have more than $t-3$ gray
neighbors in $U$ without inducing a gray book $B_{t-2}$ with $u_0$.  Therefore,

\begin{eqnarray*}
   \ell\left[\frac{p-g}{p}-x\right]+\frac{1-2p}{p}d_G(u_0)-(t-3)d_G(u_0)
   & \leq & (t-1)(1-x-d_G(u_0))\\
   \ell\left[\frac{p-g}{p}-x\right]& \leq & (t-1)(1-x)-\frac{1}{p}d_G(u_0) .
\end{eqnarray*}

Recalling that by Proposition \ref{prop:xbound}, $\frac{p-g}{p}\geq x$, we use the pigeon-hole bound from Fact~\ref{ph} $\ell\geq d_G(u_0)/x$ to get
\begin{eqnarray*}
   \frac{d_G(u_0)}{x}\left[\frac{p-g}{p}-x\right]
   & \leq & (t-1)(1-x)-\frac{1}{p}d_G(u_0) \\
   d_G(u_0)\left[\frac{p-g}{p}-x\right]
   & \leq & (t-1)x(1-x)-\frac{x}{p}d_G(u_0) \\
   d_G(u_0)\left[\frac{p-g}{p}+\frac{1-p}{p}x\right]
   & \leq & (t-1)x(1-x) .
\end{eqnarray*}

By Proposition \ref{prop:dgformula},
$$ \left[\frac{p-g}{p}+\frac{1-2p}{p}x\right]\left[\frac{p-g}{p}+\frac{1-p}{p}x\right]
   \leq (t-1)x(1-x) . $$

\noindent Collecting terms yields,
$$ \left(\frac{p-g}{p}\right)^2 +\left[\left(\frac{p-g}{p}\right)\left(\frac{2-3p}{p}\right)-(t-1)\right]x +\left[\left(\frac{1-2p}{p}\right)\left(\frac{1-p}{p}\right)+(t-1)\right]x^2\leq
0 , $$

\noindent and so minimizing the left-hand side of the inequality with respect to $x$, we have
\begin{eqnarray*}
   \left(\frac{p-g}{p}\right)^2 -\frac{\left[(t-1)-\left(\frac{p-g}{p}\right)\left(\frac{2-3p}{p}\right)\right]^2}{4\left[\left(\frac{1-2p}{p}\right)\left(\frac{1-p}{p}\right)+(t-1)\right]} & \leq & 0 \\  
   \left(\frac{p-g}{p}\right)^2\left(4t-5\right) +2\left(\frac{p-g}{p}\right)(t-1)\left(\frac{2-3p}{p}\right)-(t-1)^2 & \leq & 0 .
\end{eqnarray*}

Using the quadratic formula,
\begin{eqnarray*}
   \frac{p-g}{p} & \leq & \frac{-2(t-1)\left(\frac{2-3p}{p}\right)+\sqrt{4(t-1)^2\left(\frac{2-3p}{p}\right)^2+4(t-1)^2(4t-5)}}{2(4t-5)} \\
   p-g & \leq & \frac{t-1}{4t-5}\left[3p-2+\sqrt{(2-3p)^2+(4t-5)p^2}\right] \\
   g & \geq & p-\frac{t-1}{4t-5}\left[3p-2+2\sqrt{1-3p+(t+1)p^2}\right] .
\end{eqnarray*}
\end{proof}

The function in (\ref{eq:genLB}) achieves its maximum at $p=\frac{2t-1}{t^2+t}$, and that maximum is, in fact,
$\frac{1}{t+1}$. Hence $\ed_{\forb(K_{2,t})}(p)$ is at least $\frac{1}{t+1}$ at $p=\frac{2t-1}{t(t+1)}$ and is at
least $\frac{1}{t+1}$ at $p=\frac{2}{t+1}$. As a result of concavity,
$$ \ed_{\forb(K_{2,t})}(p) \geq \frac{1}{t+1}\qquad\mbox{for}\qquad p\in\left[\frac{2t-1}{t(t+1)},\frac{2}{t+1}\right] . $$

Equality holds whenever $t$ is odd because, in that case, Proposition~\ref{prop:t+1} gives that
$\ed_{\forb(K_{2,t})}(p) \leq 1/(t+1)$, so $p^*_\Hmc$ must be an interval.  This concludes the proof of Theorem
\ref{thm:oddt}.

If $t\geq 5$, then we can analyze the first and second derivatives, with respect to $p$, of
\begin{equation}\label{eq:psmallLB}
   p(1-p)-\left(p-\frac{t-1}{4t-5}\left[3p-2+2\sqrt{1-3p+(t+1)p^2}\right]\right) .
\end{equation}
The maximum difference between $p(1-p)$ and the lower bound in Theorem~\ref{thm:genLB} on the interval $[0,\frac{2}{t+1}]$ is $\frac{1}{t+1}$ and occurs when $p=\frac{2t-1}{t(t+1)}$. We can also see that (\ref{eq:psmallLB}) is bounded below by $\left(\frac{1}{2}-\frac{1}{t-1}\right)p(1-p)$.

\section{Upper bound constructions}\label{7}
That we have been able to determine the entire edit distance function for $\forb(K_{2,3})$ and $\forb(K_{2,4})$
raises the question of whether it might be possible to do something similar for $\forb(K_{2,t})$ when $t\geq5$.
That is, can we always find a few simple upper bounds that determine the entire edit distance function? In this
section we show that when $t\geq5$, the number and types of known upper bounds for the function increases
significantly, though this does not necessarily preclude the possibility that a few, yet to be discovered, CRG
constructions could determine the entire function.

In Section~\ref{7.1}, we revisit the work in Section~\ref{5.1} on strongly regular graphs.  In Section~\ref{7.2},
we give some constructions inspired by the analysis of triangle-free graphs in \cite{Brandt99}.~\\

\subsection{Results from strongly regular graph constructions}\label{7.1}~\\

Recall that a \emph{strongly regular graph} with parameters $(k,d,\lambda,\mu)$ is a $d$-regular graph on $k$
vertices such that each pair of adjacent vertices has $\lambda$ common neighbors, and each pair of nonadjacent
vertices has $\mu$ common neighbors.  Here we develop a function based on the existence of a strongly regular
graph.

Suppose that $K$ is a CRG with all vertices black and all edges white or gray that is derived from a
$(k,d,\lambda,\mu)$-strongly regular graph so that the edges of the strongly regular graph correspond to gray
edges of $K$. In such a case we recall from Section \ref{5.1} that
$$f_{S_{k,d,\lambda,\mu}}(p)=\frac{1}{k}+\left(\frac{k-d-2}{k}\right)p.$$

As is commonly known (see \cite{West}, for instance), if a strongly regular graph with parameters
$(k,d,\lambda,\mu)$ exists then it is necessary, though not sufficient, for
$$ d(d-\lambda-1)=\mu(k-d-1).$$

If we substitute $\lambda=t-3$ and $\mu=t-1$ in this equation and then solve for $k$, we find that
$$k=\frac{t-1+d(d+1)}{t-1},$$
\noindent and substituting these values into $f_{S_{k,d,\lambda,\mu}}(p)$ yields
$$f_{S_{k,d,\lambda,\mu}}(p)=\frac{t-1}{t-1+d(d+1)}+\left(1-\frac{(d+2)(t-1)}{t-1+d(d+1)}\right)p.$$
Fixing $p$ and minimizing $f_{S_{k,d,\lambda,\mu}}(p)$ with respect to $d$ gives the following expression:
\begin{equation}\label{eqn:theor}\frac{p(t-2)+2(t-1)}{4t-5}-\frac{2(t-1)}{4t-5}\sqrt{1-3p+(t+1)p^2},\end{equation}
\noindent which is equal to the lower bound from (\ref{eq:genLB}) in Theorem \ref{thm:genLB}.

Of course, in order to even have a chance of actually attaining (\ref{eqn:theor}) with a strongly regular graph
construction, both $d$ and $k=\frac{t-1+d(d+1)}{t-1}$ must be integers. This equation, however, provides
something of a best case scenario for strongly regular graphs, and if there is a CRG, $K$, derived from a
$(k,d,t-3,t-1)$-strongly regular graph that realizes equation (\ref{eqn:theor}), then $f_K(p)$ is tangent to the
lower bound in (\ref{eq:genLB}) at
$$ p=\frac{2d+1}{(d+1)(d+3)-t} , $$
\noindent determining the value of $\ed_{\forb(K_{2,t})}(p)$ exactly.

The remaining upper bounds in Theorem \ref{thm:sum} are the result of checking constructions from the known
strongly regular graphs listed at \cite{Brouwer}. Figure \ref{fig:chart} (see Appendix A) is a chart of the
relevant parameters and $f_K(p)$ functions for $5\leq t\leq 8$.

To our knowledge, it is not known whether, for fixed $t$, there are a finite or infinite number of
$(k,d,t-3,t-1)$-strongly regular graphs.  See Elzinga \cite{Elzinga} for values of $\lambda$ and $\mu$ for which
the number of strongly regular graphs with parameters $(k,d,\lambda,\mu)$ is known to be finite or infinite.

There is an additional construction defining the upper bound for $t=8$ in Theorem \ref{thm:cycc}, described in
the following section.

\subsection{Cycle construction}\label{7.2}
\begin{defn}[\cite{West}, p.296] For two vertices $x,y\in V(G)$, where $G$ is a simple connected graph, let
${\rm len}(x,y)$ denote the length of the minimum path from $x$ to $y$.  The $r^{\rm th}$ power of $G$, $G^r$, is
the graph with vertex set $V(G^r)=V(G)$, and edge set $E(G^r)=\{xy:x\neq y \text{ and } {\rm len}(x,y)\leq r\}$.
\end{defn}

Let $C_k^r$ be the cycle on $k$ vertices raised to the $r$th power. Define $C_{k,r}$ to be the CRG on $k$ black
vertices with white edges corresponding to those in $C_k^r$ and gray edges corresponding to those in the
complement of $C_k^r$. Recall that ${\rm \ew}$ denotes the set of white edges for a given CRG.  Then
\begin{eqnarray*}
f_{C_{k,r}}(p) & = & \frac{1}{k^2}[(1-p)k+2p\ews] \\
& = & \frac{1}{k^2}[(1-p)k+2p(rk)] \\
& = & \left(\frac{2r-1}{k}\right)p+\frac{1}{k}.
\end{eqnarray*}

\begin{prop} $C_{5+t,2}$ forbids a $K_{2,t}$ embedding, and therefore $\ed_{\forb(K_{2,t})}(p)\leq
\frac{3p+1}{5+t}$.
\end{prop}

\begin{proof}
First, we check that $C_{5+t,2}$ does not contain a gray $\Kt$. If $u_1$ and $u_2$ are any two vertices in
$C_{5+t}^2$, then $|(N(u_1)\cup N(u_2))-\{u_1,u_2\}|\geq4$.  This inequality is justified by observing that two
vertices $u_1$ and $u_2$ that are neighbors in $C_{5+t}$ have the smallest possible number of total neighbors in
$C_{5+t}^2$, and this common neighborhood has order $4$. It then follows that $|N(u_1)\cap N(u_2)|\leq t-1$ in
the complement of $\displaystyle C_{5+t}^2$.  Thus, $C_{5+t,2}$ does not contain a gray $\Kt$.

Second, we check that  $C_{5+t,2}$ does not contain a gray $\T$. If $u_1$ and $u_2$ are any two nonadjacent
vertices in $C_{5+t}^2$, then $|(N(u_1)\cup N(u_2))-\{u_1,u_2\}|\geq6$. Therefore, by reasoning similar to above,
$|N(u_1)\cap N(u_2)|\leq {t}-3$ in the complement of $C_{5+t}^2$, implying $C_{5+t,2}$ does not contain a gray
$\T$.

Thus, by Lemma \ref{lem:forbsub}, $C_{5+t,2}$ forbids a $K_{2,t}$ embedding, and therefore
$\ed_{\forb(K_{2,t})}(p)\leq f_{C_{5+t,2}}(p)=\frac{3p+1}{5+t}$.
\end{proof}

While there are several other orders and powers of cycles that would also lead to a construction forbidding
$K_{2,t}$ embedding, none of them have a corresponding $f_K(p)$ value that beats the upper bound
$\min\{p(1-p),\frac{3p+1}{5+t},\frac{1-p}{t-1}\}$, so we restrict our interest to this one.

For $t\geq 5$, $f_{C_{5+t,2}}(p)$ is always an improvement on the bound $\min\{p(1-p),\frac{1-p}{t-1}\}$ from Theorem \ref{thm:trivbds}, though it is improved upon or made irrelevant by bounds from strongly regular graphs, for $t\leq 7$. When $t=4$, the function $f_{C_{9,2}}(p)$ is tangent to the edit distance function at $p=1/3$, where the edit distance function achieves its maximum value.~\\

\subsection{F\"{u}redi constructions}~\\

As is observed in Lemma \ref{lem:forbsub} and used in the exploration of the past two constructions, graphs that
forbid $K_{2,t}$ and $B_{t-2}$ as subgraphs are of interest when looking for CRGs that forbid $K_{2,t}$
embedding.  The following results come from examining the bipartite versions of $K_{2,t}$-free graph
constructions described by F\"{u}redi \cite{FurediZar}. This strategy mimics the one used in \cite{MT} with
Brown's $K_{3,3}$-free construction.

\begin{proof}[Proof of Theorem \ref{thm:tg9}.]
We take the construction described in \cite{FurediZar} for a $K_{2,t}$-free graph $G$ on $n=(q^2-1)/(t-1)$
vertices, each with degree $q$, where $q$ is a prime power so that $t-1$ divides $q-1$.  We should note here that
in the original construction from \cite{FurediZar}, loops were omitted, reducing the degree of some vertices to
$q-1$. It is to our advantage, however, to leave the loops in so that the final construction will be $q$-regular.
By the same proof as in \cite{FurediZar}, the graph with loops still retains the property that no two vertices
have a common neighborhood greater than $t-1$ even when a looped vertex is considered to be in its own
neighborhood.

Next, we create a CRG, $K$, by taking two copies of the vertex set $\{v_1,\ldots,v_n\}$ from the $K_{2,t}$-free
graph with loops described above: $\{v_1',\ldots,v_n'\}$, $\{v_1'',\ldots,v_n''\}$. Color all of these $k=2n$
vertices black, and let $\egk=\{v_i'v_j'':v_iv_j\in E(G)\}$ with all edges not in $\egk$ white.

The gray subgraph of $K$ is bipartite, so it cannot contain a $B_{t-2}$, and since no two vertices $v_i$ and
$v_j$ from the original construction have more than $t-1$ common neighbors, the common neighborhood of two
vertices in the gray subgraph of $K$ is also at most $t-1$. Thus by Lemma \ref{lem:forbsub}, $K$ forbids a
$K_{2,t}$ embedding.

The CRG, $K$, has $k=2n=2(q^2-1)/(t-1)$ vertices and $q(q^2-1)/(t-1)$ gray edges, so by equation (\ref{eqn:f}),
$f_K(p)$ is as described in the statement of Theorem \ref{thm:tg9}.
\end{proof}

\begin{rem}
Though the property of being bipartite is sufficient to exclude a $B_{t-2}$ subgraph, using a bipartite
$K_{2,t}$-free construction may not be the optimal choice. A more efficient CRG may be constructed from another
graph that has a gray subgraph that is both $K_{2,t}$- and $B_{t-2}$-free, but, for instance, still contains
triangles.
\end{rem}

Nevertheless, we can discover more about the potential for these constructions to improve upon the bounds for
$\ed_{\forb(K_{2,t})}(p)$ by fixing $p$ and considering the general formula in Theorem \ref{thm:tg9} as a
continuous function with respect to $q$.

\begin{lem}\label{lem:q_0}
Let $t\geq 3$, and let $q_0<q$ be prime powers such that $t-1$ divides both $q_0-1$ and $q-1$.  If the CRG,
$K_0$, is constructed according to the proof of Theorem \ref{thm:tg9} with parameter $q_0$ and if the CRG, $K$,
is constructed according to the proof of Theorem \ref{thm:tg9} with parameter $q$, then $f_{K_0}(p)\leq f_{K}(p)$
for $p\in\left[\frac{2}{4+q_0},\frac{1}{3}\right)$.
\end{lem}

\begin{proof}
We begin the proof by fixing $p$ and $t$ and analyzing $\phi(q)=\frac{t-1+p(2q^2-q(t-1)-2t)}{2(q^2-1)}$.  Note
that $f_{K_0}(p)=\phi(q_0)$, and $f_K(p)=\phi(q)$.  Consider when the derivative

$$\phi'(q)= \frac{(t-1)(q^2p+p+4qp-2q)}{2(q^2-1)^2}$$
\noindent is positive and, therefore, $\phi$ is increasing.  Since the greater value of $q$ that makes
$q^2p+p+4qp-2q=0$ (note that the leading term is nonnegative) occurs at $q=\frac{(1-2p)+\sqrt{(1-2p)^2-p^2}}{p}$,
it follows that $\phi'(q)\geq 0$ when $q\geq\frac{(1-2p)+\sqrt{(1-2p)^2-p^2}}{p}$.  If $p<1/3$ and
$q_0\geq\frac{2(1-2p)}{p}$, then
$$ q>q_0\geq\frac{2(1-2p)}{p}>\frac{(1-2p)+\sqrt{(1-2p)^2-p^2}}{p} . $$

Thus, $\phi'(q)\geq 0$ for $\frac{2}{4+q_0}\leq p<1/3$.  Therefore $f_{K_0}(p)\leq f_K(p)$ for $p$ in this
interval.
\end{proof}

Additionally, we can make some statements about when we can expect constructions that originate from the
$K_{2,t}$-free graphs described by F\"uredi \cite{FurediZar} to improve upon the bound $p(1-p)$ for any $q$.

\begin{lem}\label{qvalues9}
Fix $t\geq9$, and let $q$ be a prime power such that $t-1$ divides $q-1$.  Let $K$ be the CRG with parameter $q$
described in the proof of Theorem \ref{thm:tg9}, hence $f_K(p)=\frac{t-1+p(2q^2-q(t-1)-2t)}{2(q^2-1)}$. Then for
any sufficiently large prime power $q$ and corresponding $K$, there is an interval of values of $p$ on which
$f_K(p)<p(1-p)$. Moreover as $q\rightarrow\infty$ the left-hand endpoints of these open intervals approach $0$.

That is, we can find an infinite sequence of CRG constructions that improve upon the known bounds for
$\forb(K_{2,t})$ when $t\geq 9$, and the intervals on which these improvements occur get arbitrarily close to
$0$.

\end{lem}

\begin{proof}
We begin by observing that
$f_K(p)=\frac{t-1+p(2q^2-q(t-1)-2t)}{2(q^2-1)}=p-\frac{p(q(t-1)+2t-2)-(t-1)}{2(q^2-1)}$. Thus if $f_K(p)<p(1-p)$,
\begin{eqnarray}
p-\frac{p(q(t-1)+2t-2)-(t-1)}{2(q^2-1)} & < & p-p^2 \nonumber \\
2p^2(q^2-1) & < & p(q(t-1)+2(t-1))-(t-1) \nonumber \\
2p^2(q^2-1)-p(t-1)(q+2)+(t-1) & < & 0 . \label{eqn:fur}
\end{eqnarray}

The minimum value of $2p^2(q^2-1)-p(t-1)(q+2)+(t-1)$ occurs when $p=\frac{(t-1)(q+2)}{4(q^2-1)}$.  Therefore, the
inequality above is satisfied for some $q$ and $p$ values if and only if
\begin{eqnarray*}
2\left[\frac{(t-1)(q+2)}{4(q^2-1)}\right]^2(q^2-1)-\left[\frac{(t-1)(q+2)}{4(q^2-1)}\right](t-1)(q+2)+(t-1) & < & 0 \\
(t-1)\left(1-\frac{(t-1)(q+2)^2}{8(q^2-1)}\right) & < & 0 .
\end{eqnarray*}

That is, $f_K(p)$ from the constructions in \cite{FurediZar} is less than $p(1-p)$ for some value of $p$ if and
only if $1-\frac{(t-1)(q+2)^2}{8(q^2-1)}<0$. For positive $q$, it is always the case that $(q+2)^2>q^2-1$, and so
any $q$ satisfying the constraints of the original construction will improve upon the upper bound established by
$p(1-p)$ for some $p$ when $t\geq 9$. Furthermore, for a fixed prime power $q$ for which $t-1$ divides $q-1$, it
is a definite improvement for some open neighborhood around $p=\frac{(t-1)(q+2)}{4(q^2-1)}$.  This value
approaches $0$ as $q\rightarrow\infty$, and there are an infinite number of prime powers $q$ such that $t-1$
divides $q-1$ (see \cite{FurediZar}).  Thus, it is the case that for arbitrarily small $p$, we can find some $q$
such that $f_K(p)<p(1-p)$.

\end{proof}

Lemma~\ref{qvalues8} is an analysis of the F\"uredi constructions when $t\leq 8$. We then show that our bounds
from these constructions do not have an effect on the value of $\ed_{\forb(K_{2,t})}(p)$ for $t\leq 8$.

\begin{lem}\label{qvalues8}
Fix $5\leq t\leq8$, and let $q$ be a prime power such that $t-1$ divides $q-1$.  Let $K$ be the CRG with
parameter $q$ described in the proof of Theorem \ref{thm:tg9}, hence
$f_K(p)=\frac{t-1+p(2q^2-q(t-1)-2t)}{2(q^2-1)}$. Then
\begin{equation}\label{eq:qbound}q<\frac{(t-1)+\sqrt{(t-1)^2+(9-t)(t+1)}}{\frac{1}{2}(9-t)} .\end{equation}
\end{lem}

\begin{proof} Returning to inequality (\ref{eqn:fur}) and performing a similar analysis to that in the proof of Lemma \ref{qvalues9}, we see that if
$t\leq8$, then $2p^2(q^2-1)-p(t-1)(q+2)+(t-1)<0$ for some value of $p$ if and only if
$$\frac{(t-1)-\sqrt{(t-1)^2+(9-t)(t+1)}}{\frac{1}{2}(9-t)}<q<\frac{(t-1)+\sqrt{(t-1)^2+(9-t)(t+1)}}{\frac{1}{2}(9-t)}.$$
 \noindent The lower bound for $q$ described above is immaterial since for $t\leq8$ it is always negative. The upper
 bound completes the proof of Lemma \ref{qvalues8}.
\end{proof}

Using Lemma \ref{qvalues8}, we generated the following table of possible $q$ values that obey the inequality in
(\ref{eq:qbound}).  Since we have already determined the entire edit distance function for $\forb(K_{2,3})$ and
$\forb(K_{2,4})$, only $t=5,6,7,8$ needed to be considered:
\begin{center}
\begin{tabular}{|c||c|}
\hline
  t & \text{possible $q$ values}\\ \hline
  5 & 5 \\
  6 & \text{none} \\
  7 & 7,\ 13 \\
  8 & 8,\ 29 \\
  \hline
\end{tabular}

\end{center}

A case analysis of the $f_K(p)$ functions corresponding to these $q$ values finds no improvement to the bounds
established by $\min\{p(1-p),\frac{3p+1}{t+5},\frac{1-p}{t-1}\}$, except in the cases when $t=7$ and $q=13$, and
$t=8$ and $q=29$. In these cases, we see an improvement for the approximate ranges $p\in(0.125,0.1358)$ and
$p\in(0.0625,0.06667)$, respectively, but even these improvements are surpassed by results from strongly regular
graph constructions.

\section{Conclusion}
\begin{itemize}

\item Although we determine all of $\ed_{\forb(K_{2,4})}(p)$, convexity allows $d^*_{\forb(K_{2,4})}$ to be
determined with only Lemma \ref{p large}.  Furthermore, the generalized quadrangle ${\rm GQ}(2,2)$ was
unnecessary to compute this quantity, since $p(1-p)<2/9$ for $p\in [0,1/3)$, and $p(1-p)$ is an upper bound for
every function $\ed_{\forb(K_{2,t})}(p)$, $t\geq 2$.

\item Proposition \ref{p small 2} gives a nontrivial lower bound for a black-vertex, $p$-core CRG, $K$, that
forbids a $K_{2,4}$ embedding.  If we take inequality (\ref{eq:remark}) and solve for $g$, then we see that,
if $p\in (0,1/2)$, then
$$ g_K(p)\geq\frac{2p+6-6\sqrt{1-3p+5p^2}}{11},$$
which is strictly larger than $p(1-p)$ for $p\in (0,(9-4\sqrt{3})/11)$ and corresponds to the general lower bound
described in Theorem \ref{thm:genLB}. In particular, there is a positive gap between the $g_K(p)$ functions for
black-vertex CRGs and the CRG with one white and one black vertex.

\item Ed Marchant reports having also proven that $\ed_{\forb(K_{2,3})}(p)=\min\{p(1-p),(1-p)/2\}$, using different methods.

\item While we have yet to determine the entire edit distance function $\ed_{\forb(K_{2,t})}(p)$ when
$t\geq 5$, strongly regular graph constructions have the potential to determine its value exactly, at
least for certain values of $p$.  It is likely that with the development of knowledge of strongly regular graphs,
we will see improved upper bounds for $\ed_{K_{2,t}}(p)$.  Already, they provide significant improvements to
previously known upper bounds, realizing the function value exactly in some cases.~\\

\item For $t\geq 9$, F\"{u}redi's $K_{2,t}$-free construction leads to improvements to the bound $p(1-p)$ for
values of $p$ arbitrarily close to $0$.  In fact, analyzing inequality (\ref{eqn:fur}) with respect to $q$,
indicates that for fixed $p$ and $t\geq 9$, if an appropriate $q$ exists such that
$$ \frac{t-1-\sqrt{(t-1)^2-8(t-1)(1-2p)-16p^2}}{4p}<q<\frac{t-1+\sqrt{(t-1)^2-8(t-1)(1-2p)-16p^2}}{4p} , $$
then there is an improvement.

Meanwhile, for $t\leq 8$, the upper bounds from these constructions, at least when the corresponding $f_K$
functions are considered, are inferior to those from alternative constructions.  It is unknown, whether or not
these improvements are best possible, or if there is another construction that could render them irrelevant too.~\\

\end{itemize}

\section{Acknowledgements}
\label{sec:ack} Special thanks to Sung-Yell Song for his knowledge and expertise on strongly regular graphs and
Maria Axenovich for providing references on F\"uredi's constructions. Thank you to Ed Marchant and Andrew
Thomason for helpful conversations and to Stephen Hartke for providing some references on strongly regular
graphs. We also wish to acknowledge some helpful comments from anonymous referees.


\pagebreak
\appendix
\FloatBarrier
\section{Figures}
\label{sec:figs}
\begin{figure}[htbp]
$$\begin{array}{ccc}
      \begin{tabular}{|l|c|c|}
  \hline
  $t$ values & \text{parameters} & $f_K(p)$ \\
  \hline
  $t\geq 5$ & $(13,6,2,3)$ & $(1+5p)/13$ \\
  & $(40,12,2,4)$ & $(1+26p)/40$ \\
  & $(96,19,2,4)$ & $(1+75p)/96$ \\
  \hline
  $t\geq 6$ & $(10,6,3,4)$ & $(1+2p)/10$ \\
  & $(17,8,3,4)$ &  $(1+7p)/17$\\
  & $(26,10,3,4)$ & $(1+14p)/26$ \\
  & $(85,20,3,5)$ &  $(1+63p)/85$\\
\hline
\end{tabular} &\ \ \ \ \ \ \ &\begin{tabular}{|l|c|c|}
\hline
  $t$ values & \text{parameters} & $f_K(p)$ \\
  \hline
  $t\geq 7$ & $(16,9,4,6)$ & $(1+5p)/16$ \\
  & $(36,14,4,6)$ & $(1+20p)/36$ \\
  & $(49,16,3,6)$ &  $(1+31p)/49$\\
  & $(64,18,2,6)$ &  $(1+44p)/64$\\
  & $(100,22,0,6)$ & $(1+76p)/100$ \\
  & $(156,30,4,6)$ & $(1+124p)/156$ \\
  \hline
  $t\geq 8$ & $(25,12,5,6)$ &  $(1+11p)/25$\\
  & $(76,21,2,7)$ &  $(1+53p)/76$\\
  & $(125,28,3,7)$ & $(1+95p)/125$ \\
  \hline
\end{tabular} \\
  \end{array}$$
\caption{Above are the known parameters (see \cite{Brouwer}) and $f_K(p)$ functions from strongly regular graphs
that provide an improvement upon the known upper bound for $\ed_\Hmc(p)$ for some interval of $p$ values, where
$\Hmc=\forb(K_{2,t})$ for $5\leq t\leq 8$.  Parameters with resulting bounds surpassed by other strongly regular
graph constructions are omitted.}\label{fig:chart}
  \end{figure}

  \begin{figure}[htbp]\ \hfill%
\begin{minipage}[t]{2.5in}
\includegraphics[width=2.5in]{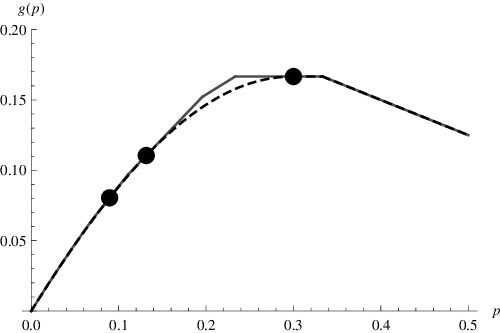}
\caption{Upper and lower bounds (in solid and dashed respectively) for $\ed_{\forb(K_{2,5})}(p)$. Points indicate
tangency.}
\end{minipage}\ \hfill\ %
\begin{minipage}[t]{2.5in}
\includegraphics[width=2.45in]{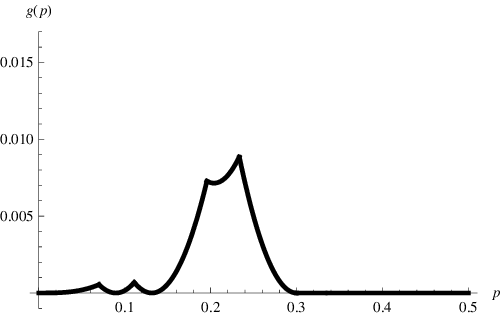}
\caption{Difference between upper and lower bounds for $\ed_{\forb(K_{2,5})}(p)$.}
\end{minipage}\ \hfill\ %
\end{figure}

\begin{figure}[htbp]\ \hfill%
\begin{minipage}[t]{2.5in}
\includegraphics[width=2.45in]{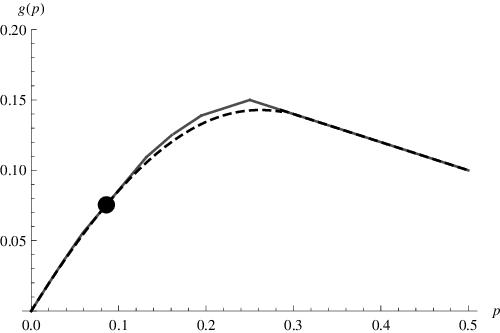}
\caption{Upper and lower bounds (in solid and dashed respectively) for $\ed_{\forb(K_{2,6})}(p)$. Points indicate
tangency.}
\end{minipage}\ \hfill\ %
\begin{minipage}[t]{2.45in}
\includegraphics[width=2.5in]{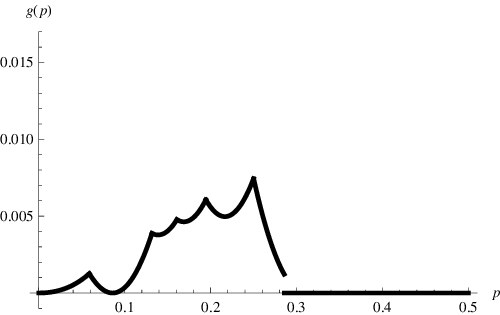}
\caption{Difference between upper and lower bounds for $\ed_{\forb(K_{2,6})}(p)$.}
\end{minipage}\ \hfill\ %
\end{figure}

\begin{figure}[htbp]\ \hfill%
\begin{minipage}[t]{2.45in}
\includegraphics[width=2.45in]{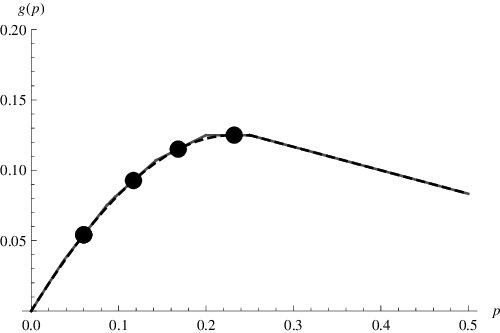}
\caption{Upper and lower bounds (in solid and dashed respectively) for $\ed_{\forb(K_{2,7})}(p)$. Points indicate
tangency.}
\end{minipage}\ \hfill\ %
\begin{minipage}[t]{2.45in}
\includegraphics[width=2.5in]{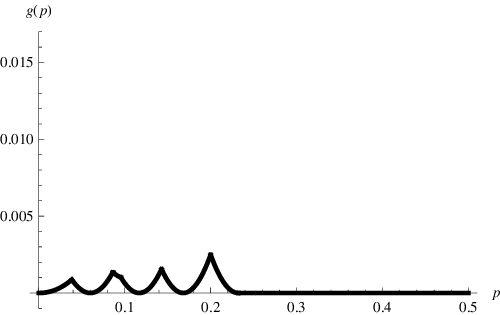}
\caption{Difference between upper and lower bounds for $\ed_{\forb(K_{2,7})}(p)$.}
\end{minipage}\ \hfill\ %
\end{figure}

\begin{figure}[htbp]\ \hfill%
\begin{minipage}[t]{2.45in}
\includegraphics[width=2.45in]{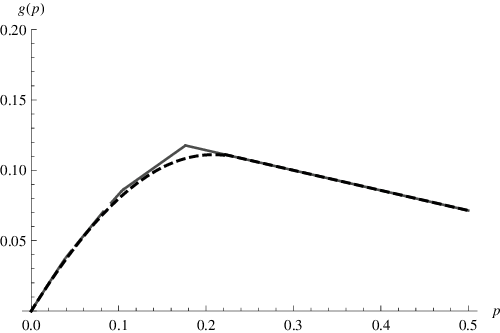}
\caption{Upper and lower bounds (in solid and dashed respectively) for $\ed_{\forb(K_{2,8})}(p)$. In this
instance there are no points of tangency.}
\end{minipage}\ \hfill\ %
\begin{minipage}[t]{2.5in}
\includegraphics[width=2.45in]{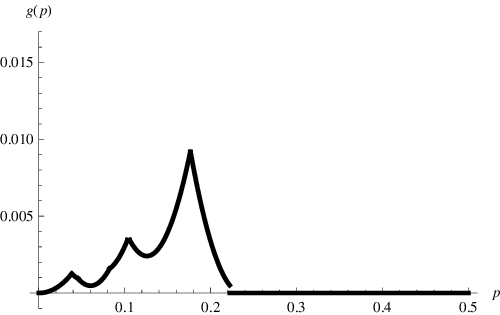}
\caption{Difference between upper and lower bounds for $\ed_{\forb(K_{2,8})}(p)$.}
\end{minipage}\ \hfill\ %
\end{figure}
\FloatBarrier
\section{Proofs of Lemma \ref{p large}, Proposition \ref{p small 1}, Proposition \ref{p small 2}, and Lemma \ref{p small}}
\label{sec:lemproofs}

\noindent \emph{\textbf{Proof of Lemma \ref{p large}:}}
\begin{proof}
We break this into two cases: when $K$ does and does not have a gray triangle. \\

\noindent\textbf{Case 1:} $K$ has a gray triangle. \\

Let the gray subgraph of $K$ contain a triangle whose vertices are $v_1$, $v_2$ and $v_3$ with optimal weights
$y_1,\ y_2$ and $y_3$, respectively. Because $K$ has no gray $B_2$, we know that no pair of the vertices
$v_1,v_2,v_3$ have a common gray neighbor other than the remaining vertex in the triangle. Letting $g=g_K(p)$, we
have the following because the sum of the optimal weights on all vertices in $K$ is $1$:

\begin{eqnarray*}
y_1+y_2+y_3+\displaystyle\sum_{i=1}^3\left[d_G(v_i)-(y_1+y_2+y_3-y_i)\right] & \leq & 1 .
\end{eqnarray*}
Then, applying Proposition \ref{prop:dgformula},
\begin{eqnarray*} y_1+y_2+y_3+3\left(\frac{p-g}{p}\right)+\frac{1-2p}{p}(y_1+y_2+y_3)-2(y_1+y_2+y_3) & \leq & 1 \\
3\left(\frac{p-g}{p}\right)+\frac{1-3p}{p}(y_1+y_2+y_3) & \leq & 1 ,
\end{eqnarray*}
and so
\begin{eqnarray*}
\frac{2p-3g}{p} & \leq & \left(\frac{3p-1}{p}\right)(y_1+y_2+y_3) \leq \frac{3p-1}{p} .
\end{eqnarray*}

Consequently, $g\geq (1-p)/3$ with equality if and only if $y_1+y_2+y_3=1$; i.e., $K$ itself is a gray triangle. \\

\noindent\textbf{Case 2:} $K$ has no gray triangle. \\

Let $u_0$ be a vertex of largest weight, $x=\xb(u_0)$, and let $U=N_G(u_0)$.  The absence of a gray triangle
means that there are no gray edges between pairs of vertices in $U$. Furthermore, no vertex in $W$ can be
adjacent to more than three vertices in $U$ via a gray edge, since by Lemma \ref{lem:forbsub}, the gray subgraph
of $K$ does not contain a $K_{2,4}$.

Let $u_1,\ldots,u_\ell$ be an enumeration of the vertices in $U$ with weights $x_1,\ldots,x_\ell$, respectively,
and $g=g_K(p)$. Then
\begin{eqnarray*}
\sum_{i=1}^\ell \left(d_G(u_i)-x\right) & \leq & 3\xb(W) \\
& \leq & 3(1-x-\xb(U)),
\end{eqnarray*}
and applying Proposition \ref{prop:dgformula} to compute $d_G(u_i)$,
\begin{eqnarray}
\sum_{i=1}^\ell\left(\frac{p-g}{p}+\frac{1-2p}{p}x_i-x\right) & \leq & 3(1-x-\xb(U)) \nonumber \\
\ell\left(\frac{p-g}{p}-x\right)+\frac{1-2p}{p}\xb(U) & \leq & 3(1-x)-3\xb(U) \nonumber \\
\ell\left(\frac{p-g}{p}-x\right) & \leq & 3(1-x)-\frac{1+p}{p}\xb(U) . \label{eqn:3}
\end{eqnarray}

First, suppose $\ell\geq 5$.  Then, from inequality (\ref{eqn:3}), we have
\begin{eqnarray*}
5\left(\frac{p-g}{p}-x\right) & \leq & 3(1-x)-\frac{1+p}{p}\xb(U),
\end{eqnarray*}
\noindent and applying Proposition \ref{prop:dgformula} again,

\begin{eqnarray*}
5\left(\frac{p-g}{p}-x\right) & \leq & 3(1-x)-\frac{1+p}{p}\left(\frac{p-g}{p}+\frac{1-2p}{p} x\right)\\
\frac{1+6p}{p}\cdot\frac{p-g}{p}-3 & \leq & \left(5-3-\frac{1+p}{p}\cdot\frac{1-2p}{p}\right)x \\
p(1+3p)-g(1+6p) & \leq & x\left(4p^2+p-1\right) .
\end{eqnarray*}

If $4p^2+p-1<0$, then we may use the fact that $x>0$,
$$ g>\frac{p(1+3p)}{1+6p}= \frac{1-p}{3}+\frac{(3p-1)(1+5p)}{3(1+6p)}. $$
If $4p^2+p-1\geq 0$, then we use Proposition \ref{prop:xbound} and substitute $x=g/(1-p)$,
\begin{eqnarray*}
p(1+3p)-g(1+6p) & \leq & \frac{g}{1-p}\left(4p^2+p-1\right) \\
p(1+3p) & \leq & g\left(\frac{6p-2p^2}{1-p}\right) \\
\frac{1-p}{3}+\frac{(1-p)(11p-3)}{6(3-p)} & \leq & g .
\end{eqnarray*}
Regardless of the value of $p\in (1/3,1/2)$, if $\ell\geq 5$, then $g>(1-p)/3$. Therefore, we may assume that
$\ell\leq 4$.

Second, suppose $\ell\leq 2$.  Then by Fact \ref{ph} we have $\ell\geq\xb(U)/x$, yielding

$$ \xb(U)/x \leq \ell \leq  2,$$
\noindent and so bounding $\xb(U)$ using Proposition \ref{prop:dgformula},
\begin{eqnarray*}
\frac{1}{x}\left(\frac{p-g}{p}+\frac{1-2p}{p}x\right)&\leq& 2\\
\frac{p-g}{p} & \leq & \frac{4p-1}{p}x.
\end{eqnarray*}
\noindent Using Proposition \ref{prop:xbound}, $x\leq \frac{g}{1-p}$ yields
\begin{eqnarray*}
\frac{p-g}{p} & \leq & \frac{4p-1}{p}\cdot\frac{g}{1-p}\\
p(1-p)&\leq&3pg,
\end{eqnarray*}
and so if $\ell\leq 2$, then $g\geq (1-p)/3$, with equality if and only if $x=g/(1-p)$, and consequently, $K$ is
a gray triangle. So, we may further assume that $\ell\in\{3,4\}$.

Third, suppose $\ell=3$.  Then
\begin{eqnarray*}
\xb(U)/x & \leq & 3 \\
\frac{p-g}{p} & \leq & \frac{5p-1}{p}x \\
\frac{p-g}{5p-1} & \leq & x .
\end{eqnarray*}
Returning to inequality (\ref{eqn:3}), we have
\begin{eqnarray*}
3\left(\frac{p-g}{p}-x\right) & \leq & 3(1-x)-\frac{1+p}{p}\xb(U) \\
\frac{1+4p}{p}\cdot\frac{p-g}{p}-3 & \leq & -\left[\frac{1+p}{p}\cdot\frac{1-2p}{p}\right]x \\
p(1+p)-g(1+4p) & \leq & -(1+p)(1-2p)\left(\frac{p-g}{5p-1}\right) \\
p(1+p)(5p-1)+p(1+p)(1-2p) & \leq & g\left[(1+p)(1-2p)+(1+4p)(5p-1)\right] \\
\frac{1+p}{6} & \leq & g \\
\frac{1-p}{3}+\frac{3p-1}{6} & \leq & g .
\end{eqnarray*}
If $\ell=3$, then $g>(1-p)/3$.

Fourth, and finally, suppose $\ell=4$.  Then
\begin{eqnarray*}
\xb(U)/x & \leq & 4 \\
\frac{p-g}{p} & \leq & \frac{6p-1}{p}x \\
\frac{p-g}{6p-1} & \leq & x .
\end{eqnarray*}
Returning to inequality (\ref{eqn:3}), we have
\begin{eqnarray*}
4\left(\frac{p-g}{p}-x\right) & \leq & 3(1-x)-\frac{1+p}{p}\xb(U) \\
\frac{1+5p}{p}\cdot\frac{p-g}{p}-3 & \leq & \left[4-3-\frac{1+p}{p}\cdot\frac{1-2p}{p}\right]x \\
p(1+2p)-g(1+5p) & \leq & \left[3p^2+p-1\right]x .
\end{eqnarray*}

If $3p^2+p-1<0$, then we use the fact that $x\geq (p-g)/(6p-1)$:
\begin{eqnarray*}
p(1+2p)-g(1+5p) & \leq & \left[3p^2+p-1\right]\left[\frac{p-g}{6p-1}\right] \\
p(1+2p)-\frac{p(3p^2+p-1)}{6p-1} & \leq & g\left[1+5p-\frac{3p^2+p-1}{6p-1}\right] \\
\frac{1+3p}{9} & \leq & g \\
\frac{1-p}{3}+\frac{2(3p-1)}{9} & \leq & g .
\end{eqnarray*}
If $3p^2+p-1\geq 0$, then we use Fact \ref{ph} to bound $x\leq\frac{g}{1-p}$,
\begin{eqnarray*}
p(1+2p)-g(1+5p) & \leq & \left[3p^2+p-1\right]\left[\frac{g}{1-p}\right] \\
p(1+2p) & \leq & g\left[1+5p+\frac{3p^2+p-1}{1-p}\right] \\
\frac{(1-p)(1+2p)}{5-2p} & \leq & g \\
\frac{1-p}{3}+\frac{2(4p-1)(1-p)}{3(5-2p)} & \leq & g .
\end{eqnarray*}
Regardless of the value of $p\in (1/3,1/2)$, if $\ell=4$, then $g>(1-p)/3$.

This ends Case 2 and the proof of the lemma.\\

\end{proof}

\noindent \emph{\textbf{Proof of Proposition \ref{p small 1}:}}
\begin{proof}
Let $u_1,\ldots,u_\ell$ be an enumeration of the vertices of $U$.  Then
\begin{eqnarray*}
   \sum_{i=1}^{\ell}\left(d_G(u_i)-x\right)-\xb(U_1) & \leq & 3(1-x-\xb(U)),
\end{eqnarray*}
   \noindent and applying Proposition \ref{prop:dgformula},
\begin{eqnarray*}
   \sum_{i=1}^{\ell}\left(\frac{p-g}{p}+\frac{1-2p}{p}\xb(u_i)-x\right)-\xb(U_1) & \leq & 3(1-x-\xb(U)) .
\end{eqnarray*}
Simplification yields the first inequality.  The second inequality results from observing that
$\xb(U_1)\leq\xb(U)$.\\

\end{proof}

\noindent \emph{\textbf{Proof of Proposition \ref{p small 2}:}}
\begin{proof}
We begin with Proposition \ref{p small 1} and then use $\ell\geq\xb(U)/x$ from Fact \ref{ph}:
\begin{eqnarray}
   \ell\left(\frac{p-g}{p}-x\right) & \leq & 3-3x-\frac{1}{p}\xb(U) \nonumber \\
   \frac{\xb(U)}{x}\left(\frac{p-g}{p}-x\right) & \leq & 3-3x-\frac{1}{p}\xb(U) \nonumber \\
   \xb(U)\left(\frac{p-g}{px}-1+\frac{1}{p}\right) & \leq & 3-3x \nonumber \\
   \left[\frac{p-g}{p}+\frac{1-2p}{p}x\right] \left[\frac{p-g}{p}+\frac{1-p}{p}x\right] & \leq & 3x-3x^2. \label{eq:remark}
   \end{eqnarray}
\noindent  Recall that $(p-g)/p\geq p$ because $g\leq p(1-p)$, so
   \begin{eqnarray}
   \left[p+\frac{1-2p}{p}x\right] \left[p+\frac{1-p}{p}x\right] & \leq & 3x-3x^2 \nonumber \\
   p^2-(1+3p)x+\frac{1-3p+5p^2}{p^2}x^2 & \leq & 0 . \nonumber
\end{eqnarray}

The quadratic formula gives that not only must the discriminant be nonnegative (requiring $p\geq
(9-4\sqrt{3})/11$), but also
$$ x\geq\frac{p^2}{2(1-3p+5p^2)}\left[1+3p-\sqrt{-3+18p-11p^2}\right] . $$
Some routine but tedious calculations demonstrate that, for $p\in\left[(9-4\sqrt{3})/11,1/2\right)$, this expression is at least $1/25$, achieving
that value uniquely at $p=1/5$.\\

\end{proof}
\newpage
\noindent \emph{\textbf{Proof of Lemma \ref{p small}:}}
\begin{proof}
We assume that $g_K(p)\leq p(1-p)$.\\

\noindent\textbf{Case 1:} $\ell\geq 8$. \\

According to Proposition \ref{p small 1},
\begin{eqnarray*}
8\left(\frac{p-g}{p}-x\right) & \leq & \ell\left(\frac{p-g}{p}-x\right)\leq 3-3x-\frac{1}{p}\left(\frac{p-g}{p}+\frac{1-2p}{p}x\right) \\
(1-2p-5p^2)x & \leq & 3p^2-(p-g)(1+8p),
\end{eqnarray*}
\noindent and since $x\geq 1/25$ and $p-g\geq p^2$,
\begin{eqnarray*}
\frac{1-2p-5p^2}{25} & \leq & 3p^2-p^2(1+8p) \\
(1-5p)^2(1+8p) & \leq & 0,
\end{eqnarray*}
\noindent a contradiction. So, $\ell<8$. \\

\noindent\textbf{Case 2:} $\ell\leq 7$ and $x<p^2/(9p-1)$. \\

Using Fact \ref{ph}, and then Proposition \ref{prop:dgformula}
$$ 7 \geq \ell \geq \frac{\xb(U)}{x} \geq \frac{p}{x}+\frac{1-2p}{p} > \frac{9p-1}{p}+\frac{1-2p}{p} = 7 , $$
a contradiction. \\

\noindent\textbf{Case 3:} $\ell\leq 7$ and $p^2/(9p-1)\leq x\leq p/3$. \\

First we bound $\ell$:
$$ \ell\geq\frac{\xb(U)}{x}\geq\frac{p}{x}+\frac{1-2p}{p}\geq 3+\frac{1}{p}-2>6 . $$
So, $\ell=7$.  Since $\ell$ is odd, $\xb(U_1)\leq 6x$.  By Proposition \ref{p small 1},
\begin{eqnarray*}
   \ell\left(\frac{p-g}{p}-x\right) & \leq & 3-3x-\frac{1+p}{p}\xb(U)+\xb(U_1),
\end{eqnarray*}
\noindent and applying Proposition \ref{prop:dgformula},
\begin{eqnarray*}
   7\left(\frac{p-g}{p}-x\right) & \leq & 3-3x-\frac{1+p}{p}\left[\frac{p-g}{p}+\frac{1-2p}{p}x\right]+6x \\
   \frac{1-p-12p^2}{p^2}x & \leq & 3-\frac{1+8p}{p}\cdot\frac{p-g}{p} \\
   \frac{1-p-12p^2}{p^2}\left[\frac{p^2}{9p-1}\right] & \leq & 3-\frac{1+8p}{p}\cdot p \\
   \frac{(1-4p)(1+3p)}{9p-1} & \leq & 2(1-4p) \\
   \frac{1+3p}{9p-1} & \leq & 2 ,
\end{eqnarray*}
which implies $p\geq 1/5$, a contradiction. \\

\noindent\textbf{Case 4:} $\ell\leq 7$ and $x>p/3$. \\

Now we compute a stronger bound on $U_1$.  Let $u_1$ and $u_2$ be vertices in $U$ that are adjacent via a gray
edge, and let their weights be $x_1$ and $x_2$, respectively.  Then
\begin{eqnarray*}
   x+\xb(U)+(d_G(u_1)-x-x_2)+(d_G(u_2)-x-x_1) & \leq & 1
\end{eqnarray*}
because $u_1$ and $u_2$ have no common gray neighbor other than $u_0$ and because they can have no additional gray neighbor in $U$.  Applying Proposition \ref{prop:dgformula},
\begin{eqnarray*}
   x+\frac{p-g}{p}+\frac{1-2p}{p}x+2\frac{p-g}{p}-2x+\frac{1-3p}{p}(x_1+x_2) & \leq & 1 \\
   \frac{1-3p}{p}(x_1+x_2) & \leq & \frac{3g-2p}{p}-\frac{1-3p}{p}x,
   \end{eqnarray*}
   \noindent and since $p(1-p)\geq g$,
   \begin{eqnarray*}
   x_1+x_2 & \leq & p-x .
\end{eqnarray*}

We can bound the number of vertices in $U-U_1$ by using the fact that $(\ell-\ell_1)x\geq\xb(U)-\xb(U_1)$.
Returning to Proposition \ref{p small 1},
\begin{eqnarray*}
   \ell\left(\frac{p-g}{p}-x\right) & \leq & 3-3x-\frac{1+p}{p}\xb(U)+\xb(U_1) \\
   \left[\ell_1+\frac{1}{x}\xb(U)-\frac{1}{x}\xb(U_1)\right]\left(\frac{p-g}{p}-x\right)
   & \leq & 3-3x-\frac{1+p}{p}\xb(U)+\xb(U_1) \\
   \xb(U)\left(\frac{p-g}{px}-1+\frac{1+p}{p}\right)-3+3x & \leq & \xb(U_1)\left(1+\frac{p-g}{px}-1\right)-\ell_1\left(\frac{p-g}{p}-x\right).
\end{eqnarray*}

If $\ell_1=|U_1|$, then $\xb(U_1)\leq (\ell_1/2)(p-x)$.  Of course, $\xb(U)$ is lower-bounded by Proposition~\ref{prop:dgformula}.
\begin{eqnarray*}
   \left[\frac{p-g}{p}+\frac{1-2p}{p}x\right]\left(\frac{p-g}{px}+\frac{1}{p}\right)-3+3x & \leq & \frac{\ell_1}{2}(p-x)\left(\frac{p-g}{px}\right)-\ell_1\left(\frac{p-g}{p}-x\right) \\
   \left[\frac{p-g}{p}+\frac{1-2p}{p}x\right] \left(\frac{p-g}{px}+\frac{1}{p}\right)-3+3x & \leq & \ell_1\left[x-\frac{p-g}{p}\cdot\frac{3x-p}{2x}\right] \\
   \left[p+\frac{1-2p}{p}x\right]\left(\frac{p}{x}+\frac{1}{p}\right)-3+3x & \leq & \ell_1\left[x-\frac{p(3x-p)}{2x}\right] \\
   p^2-(1+2p)x+\frac{1-2p+3p^2}{p^2}x^2 & \leq & \ell_1\frac{(p-x)(p-2x)}{2} . \\
\end{eqnarray*}
Now, we bound $\ell_1$, depending on the sign of $p-2x$, requiring two more cases. \\

\noindent\textbf{Case 4a:} $\ell\leq 7$ and $x>p/3$ and $p-2x\geq 0$. \\

Here we use the bound $\ell_1\leq 6$:
\begin{eqnarray*}
   p^2-(1+2p)x+\frac{1-2p+3p^2}{p^2}x^2 & \leq & 3(p-x)(p-2x) \\
   -2p^2+(7p-1)x+\frac{1-2p-3p^2}{p^2}x^2 & \leq & 0 .
\end{eqnarray*}
By Proposition \ref{p small 2}, we may restrict our attention to $p\geq (9-4\sqrt{3})/11>1/7$ and so we may
substitute the smallest possible value for $x$, which still maintains the inequality.
\begin{eqnarray*}
   -2p^2+(7p-1)\left(\frac{p}{3}\right) +\frac{1-2p-3p^2}{p^2}\left(\frac{p}{3}\right)^2 & < & 0 \\
   -18p^2+3(7p-1)p+(1-2p-3p^2) & < & 0 \\
   1-5p & < & 0 ,
\end{eqnarray*}
a contradiction. \\

\noindent\textbf{Case 4b:} $\ell\leq 7$ and $x>p/3$ and $p-2x<0$. \\

Here we use the bound $\ell_1\geq 0$ and then replace $x$ with $\frac{p^2(1+2p)}{2-4p+6p^2}$, the value that
minimizes the left-hand side:
\begin{eqnarray*}
   p^2-(1+2p)x+\frac{1-2p+3p^2}{p^2}x^2 & \leq & 0 \\
   p^2-\frac{(1+2p)^2p^2}{4(1-2p+3p^2)} & \leq & 0 \\
   \frac{p^2(3-12p+8p^2)}{4(1-2p+3p^2)} & \leq & 0 .
\end{eqnarray*}
This, too, is a contradiction for $p\in (0,1/5)$, completing the proof of Lemma \ref{p small}.
\end{proof}

\end{document}